\documentclass[final,3p,times,twocolumn]{elsarticle}

\usepackage{url}
\usepackage{amsmath,amssymb,amsthm,array,bibentry,caption,dsfont,euscript,framed,graphicx,hyphenat,mathrsfs,natbib,remreset,totcount,xspace,picinpar,algorithmic,algorithm}
\usepackage{nicefrac}

\usepackage{geometry}
\usepackage[normalem] {ulem}

\makeatletter 
\def\ps@pprintTitle{%
 \let\@oddhead\@empty
 \let\@evenhead\@empty
 \def\@oddfoot{}%
 \let\@evenfoot\@oddfoot}
\makeatother

 \usepackage{setspace}
\let\Algorithm\algorithm
\renewcommand\algorithm[1][]{\Algorithm[#1]\setstretch{1.3}}

\newlength\myindent
\newcommand\bindent{%
 \begingroup
 \setlength{\itemindent}{\myindent}
 \addtolength{\algorithmicindent}{\myindent}
}
\newcommand\eindent{\endgroup}

\newcommand{\R}{\ensuremath{\mathbb{R}}}

\newcommand{\E}{\ensuremath{\mathbb{E}}}

\DeclareMathOperator*{\Var}{Var}

\renewcommand{\vec}[1]{{\mathbf #1}}

\begin{document}

\begin{frontmatter}

\title{Powerful statistical inference for nested data using sufficient summary statistics} 

\author[tub]{Irene Dowding\corref{cor1}}
\ead{irenedowding@web.de}

\author[tub]{Stefan Haufe\corref{cor1}}
\ead{stefanhaufe@gmail.com}

\cortext[cor1]{Authors contributed equally. Corresponding author.}

\address[tub]{Technische Universit\"at Berlin}

\begin{abstract}

Hierarchically-organized data arise naturally in many psychology and neuroscience studies. As the standard assumption of independent and identically distributed samples does not hold for such data, two important problems are to accurately estimate group-level effect sizes, and to obtain powerful statistical tests against group-level null hypotheses. A common approach is to summarize subject-level data by a single quantity per subject, which is often the mean or the difference between class means, and treat these as samples in a group-level t-test. This `naive' approach is, however, suboptimal in terms of statistical power, as it ignores information about the intra-subject variance. To address this issue, we review several approaches to deal with nested data, with a focus on methods that are easy to implement. {With what we call the sufficient-summary-statistic approach, we highlight a computationally efficient technique that can improve statistical power by taking into account within-subject variances}, and we provide step-by-step instructions on how to apply this approach to a number of frequently-used measures of effect size. The properties of the reviewed approaches and the potential benefits over a group-level t-test are quantitatively assessed on simulated data and demonstrated on EEG data from a simulated-driving experiment.
\end{abstract}


\end{frontmatter}

\section{Introduction}

Data with nested (hierarchical) structure arise naturally in many fields. In psychology and neuroimaging, for example, multiple data points are often acquired for the same subject throughout the course of an experiment; thus, there exists a subject (lower) and a group (higher) level in the data hierarchy. Two important questions are how to obtain precise estimators for group-level effect sizes from nested data, and how to obtain powerful statistical tests for the presence of group-level effects. The main difficulty associated with such nested data is that the assumption of identically distributed observations is typically violated: while samples acquired from the same subject can be considered to be identically distributed, different distributions must be assumed for different subjects. Therefore, simply pooling the data of all subjects in order to apply a standard statistical test like a t-test would lead to wrong results.


A flexible way to model nested data is to combine the data of all subject in a single linear model, referred to as the nested linear model, hierarchical linear model, multi-level model or linear mixed model \citep{Quene2004,hox2010multilevel,Woltman2012,chen2013linear}. Parameter estimation in such models is, however, difficult to implement and computationally expensive, as it typically requires non-linear optimization of non-convex objective functions. Moreover, the range of effects that can be modeled is limited to linear coefficients. It is, therefore, worthwhile to study how group-level inference can be implemented for other commonly used effect size measures such as correlations or differences in the general central tendencies of distributions. 


In the neuroimaging (e.g., electro- and magnetoencephalography, EEG/MEG) literature, the use of suboptimal inference procedures is currently still widespread, as discussed in \cite{Mumford2009,Pernet2011}. Common hierarchical approaches often summarize subject-level data by a single quantity per subject, which is often the mean or the difference between class means, and treat these as single samples in a group-level test. This `naive' summary-statistics approach is, however, not optimal in terms of statistical power, as it ignores information about the intra-subject variance. Given the low signal-to-noise ratios and small sample regimes that are typical for neuroimaging studies, the potential loss of statistical power is unfortunate.

Group-level statistical power can be improved by incorporating variances at the lower level in relatively simple ways. The problem of estimating group-level effect sizes and estimating their statistical significance can, moreover, be formulated in a compellingly simple framework, where group-level inference is conducted using the \emph{sufficient} summary statistics of separate subject-level analyses. The resulting statistical methods are simple to implement, computationally efficient, and can be easily extended to settings with more than two nesting levels, which are common, e.g., in the analysis of functional magnetic resonance imaging (fMRI) data. 
Sufficient summary statistics approaches are popular in the field of meta analysis \citep{Borenstein2009,Card2011}. In neuroimaging, they are commonly used to estimate group-level coefficients of hierarchical linear model \citep[see][for methodological reviews]{Beckmann2003,Monti2011}. Here, we argue that a wider range of popular effect size measures can benefit from the high statistical power of sufficient-summary-statistic-based estimators. While this has been exploited in various experimental studies \citep{schubert2009now,haufe2011eeg,Winkler2015,Lur20162538,batista2016developmental}, the theoretical grounds on which such estimators are derived for different effect size measures have not yet been summarized in a single accessible source.


With this paper, we aim to fill this gap by providing a review of ways to estimate group-level effect sizes and to assess their statistical significance in the context of neuroimaging experiments. We first provide a reference for a number of popular parametric and non-parametric effect size measures (Section~\ref{sec:effectsizes}), which may be skipped by readers who want to proceed directly to the nested setting. We then discuss the need to choose an appropriate group-level model, as between-subject variability differs depending on whether a `random effects' or `fixed effect' model is assumed (Section~\ref{sec:group_inference}). We also demonstrate why the simple approach of ignoring the group structure by pooling the data of all subject is invalid (Section~\ref{sec:pooling}). We then outline the popular `naive' summary-statistic approach of computing effect sizes on the subject level and treating these effect sizes as single samples in a group-level test (Section \ref{sec:naive}). {With what we call the sufficient-summary-statistic approach, we then discuss a family of techniques capable to yield unbiased group-level effect size estimates and powerful statistical tests of group-level null hypotheses, and we highlight a particular approach that yields minimum-variance effect size estimates by weighting each effect with the inverse of its variance.} In a tutorial style, we outline the steps that are required to apply this approach to different effect size measures (Section~\ref{sec:anylincomb}). Lastly, we discuss the advantages and drawbacks of Stouffer's method of combining subject-level p-values (Section~\ref{sec:stouffer}) in relation to summary-statistic approaches. 

Using synthetic data representing a two-sample separation problem, we empirically assess the performance of the reviewed approaches (Section~\ref{sec:simulations}). The properties of the various approaches and the advantages of the sufficient-summary-statistic approach are further highlighted in an application to EEG data acquired during simulated emergency braking in a driving simulator (Section~\ref{sec:realData}). {All data are provided in Matlab format along with corresponding analysis code\footnote{\url{https://github.com/stefanhaufe/GroupStats}}.}

The paper ends with a discussion of nested linear models, of multivariate extensions, and a note on non-parametric (bootstrapping and surrogate data) approaches (Section~\ref{sec:discussion}).

\section{Theory}
\label{sec:theory}

\subsection{Statistical terminology}

An \emph{effect size} $\theta$ is any quantitative measure that reflects the magnitude of some phenomenon of interest (e.g., a parameter in a model). An \emph{estimator} $\hat{\theta}$ for $\theta$ is \emph{unbiased}, if its expected value is $\theta$. 

A \emph{statistical test} is a procedure to decide, based on observed data, whether a hypothesis about a population is true. In this paper, our goal is to make inference about the presence or absence of an effect in the population. The \emph{null hypothesis} is that no effect is present. The zero effect is denoted by $\theta_0$. The null hypothesis of no effect is denoted by $H_0: \theta = \theta_0$. The \emph{alternative hypothesis} that an effect is present is denoted by $H_1$. A \emph{one-tailed} alternative hypothesis assumes that either $H_1: \theta > \theta_0$ or $H_1: \theta < \theta_0$, while a \emph{two-tailed} alternative hypothesis assumes that $H_1: \theta \neq \theta_0$.

{A \emph{test statistic} needs to be derived from the observed effect size, where its distribution under the null hypothesis is known or can be reasonably well approximated. The \emph{p-value} denotes the probability of obtaining a result at least as strong as the observed test statistic under the assumption of the null hypothesis.}
Denoting the test statistic by $T$, its cumulative distribution function under the $H_0$ by $F_T$, and its observed value in a given sample by $\tau$, the p-values for a \emph{one-tailed} alternative hypothesis are given by
\begin{align}\label{eq:one-tailed}
Pr(T \leq \tau \, | \, H_0) &= F_T(\tau) \\
Pr(T \geq \tau \, | \, H_0) &= 1 - F_T(\tau) \;,
\end{align} 
where $Pr(\cdot)$ denotes probability. The p-value for a two-tailed alternative hypothesis is given by 
%
\begin{align}\label{eq:two-tailed}
2 \cdot \min \Big\{ Pr(T \leq \tau \, | \, H_0) ,  Pr(T \geq \tau \, | \, H_0) \Big\} \;.
\end{align}

The null hypothesis is rejected if the p-value falls below an \emph{alpha-level} $\alpha$. In the opposite case, no conclusion is drawn. The most commonly used alpha-levels are $\alpha=0.05$ and $\alpha=0.01$. If the null hypothesis is rejected, we speak of a \emph{statistically significant} effect. The value of the test statistic that is required for a significant effect is called \emph{critical value}. The \emph{power} of a statistical test is the probability that the test correctly rejects the null hypothesis when the alternative hypothesis is true. Conversely, the \emph{false positive rate} is the fraction of non-existent effects that are statistically significant.
%
%


\subsection{Common effect size measures}\label{sec:effectsizes}

Before introducing the nested data setting, we review a number of popular effect size measures. For each measure, we also present an analytic expression of its variance, which is a prerequisite for assessing its statistical significance. We will later need the variance for performing statistical inference in the nested setting, too. 


\subsubsection{Mean of a sample}
A common measure of effect size is the mean of a sample. Consider a neuroimaging experiment, in which the participant is repeatedly exposed to the same stimulus. A common question to ask is whether this stimulus evokes a brain response that is significantly different from a baseline value. Assume that we observe $N$ independent samples $x_{1}, \ldots, x_{N} \in \R$. The sample mean is denoted by $\bar{x} = \nicefrac{1}{N} \sum_{i=1}^{N} x_{i}$, and the unbiased sample variance is given by $\hat{\sigma}_{\!x}^2 = \nicefrac{1}{(N-1)} \sum_{i=1}^{N} (x_{i} - \bar{x})^2$. The variance of $\bar{x}$ is given by
\begin{align}\label{eq:meanDifferenceVarianceOneSample}
\widehat{\Var}(\bar{x}) = \frac{\hat{\sigma}^2_{\!x}}{N}\;.
\end{align}
Assuming \emph{independent and identically distributed} (i.i.d.) samples, which are either normal (Gaussian) distributed or large enough,  the null hypothesis $H_0: \bar{x} = \theta_0$ can be tested using that the statistic 
\begin{align}\label{eq:tOneSample}
t = \frac{\bar{x} - \theta_0}{\surd \widehat{\Var}(\bar{x}) }
\end{align}
is approximately Student-t-distributed with $N-1$ degrees of freedom. This is the \emph{one-sample t-test}.

A similar effect size is the mean difference $\overline{x-y} = \nicefrac{1}{N} \sum_{i=1}^{N} x_{i} - y_{i}$ of two paired samples $(x_{1}, y_{1}), \ldots, (x_{N}, y_{N}) \in \R^2$. Here, the $y_i$ could, for example, represent baseline activity that is measured in each repetition right before the presentation of the experimental stimulus. A natural null hypothesis is that the mean difference is zero, i.e., $H_0: \overline{x-y} = 0$. This hypothesis can be tested with a \emph{paired t-test}, which replaces $x$ by $x - y$ in  Eq.~\eqref{eq:meanDifferenceVarianceOneSample},\eqref{eq:tOneSample}. 

Note that, if $x$ or $y$ cannot be assumed to be normal distributed, a more appropriate test is the non-parametric \emph{Wilcoxon signed-ranked test} {which tests whether the mean population ranks differ \citep{Wilcoxon1945}. Alternative robust techniques can lead to a more detailed understanding of how the groups differs, see, e.g., \cite{Rousselet2017} for a recent summary.}

\subsubsection{Difference between class-conditional means} 
\label{sec:unpairedttest}
A slightly different treatment is required for the difference between the means of two unpaired samples. Consider an experiment with two conditions $\mathcal{X}$ and $\mathcal{Y}$. In neuroimaging studies, these could differ in the type of stimulus presented. We observe $N_{\mathcal X}$ samples $x_{1}, \ldots, x_{N_{\mathcal X}} \in \R$ of brain activity within condition $\mathcal{X}$, and $N_{\mathcal Y}$ samples $y_{1}, \ldots, y_{N_{\mathcal Y}} \in \R$ within condition $\mathcal{Y}$. The sample means are denoted by $\bar{x}$ and $\bar{y}$, and their difference is given by
\begin{align}\label{eq:meanDifference}
\hat{d} = \bar{x} - \bar{y} \;.
\end{align}
%
The variance of $\hat{d}$ is estimated as 
\begin{align}\label{eq:meanDifferenceVariance}
\widehat{\Var}(\hat{d}) = \frac{\hat{\sigma}^2_{\!x}}{N_{{\mathcal X}}}+ \frac{\hat{\sigma}^2_{\!y}}{N_{{\mathcal Y}}} \;,
\end{align}
where $\hat{\sigma}^2_{\!x}$ and $\hat{\sigma}^2_{\!y}$ are the unbiased sample variances of $\mathcal{X}$ and $\mathcal{Y}$.
The null hypothesis of equal means is given by $H_0: d = 0$. Under the assumption of either normal distributed $x_{i}$ and $y_{i}$, or large enough samples, the null hypothesis can be tested with \emph{Welch's two-sample t-test}. It computes the test statistic
\begin{align}
t = \frac{\hat{d}}{\surd \widehat{\Var}(\hat{d}) }\,
\end{align}
which is approximately Student-t-distributed. The degrees of freedom can be approximated using the Welch-Satterthwaite equation \citep{welch1947generalization}. Note that assuming equal variances of $\mathcal{X}$ and $\mathcal{Y}$ leads to the better known Student's t-test, which is, however, less recommendable than Welch's t-test \citep{Ruxton2006}.

\subsubsection{Area under the ROC curve}

In many cases, one may be interested in quantifying the predictive accuracy of a binary classifier to separate experimental condition $\mathcal{X}$ from condition $\mathcal{Y}$.
{A host of evaluation criteria are available for this task, and we refer the interested reader to \cite{Baldi2000} for a comprehensive review.}
The \emph{receiver operating characteristic (ROC)} is a plot that visualizes the performance of such a binary classification system. It is obtained by plotting the true positive rate (TPR) against the false positive rate (FPR) when varying the threshold that divides the predicted condition into $\mathcal{X}$ and $\mathcal{Y}$. Assume without loss of generality that condition $\mathcal{X}$ is associated with a positive label indicating that the detection of instances of that condition is of particular interest, while $\mathcal{Y}$ is associated with a negative label. TPR is defined as the fraction of correctly classified positive samples among all positive samples, while FPR denotes the fraction negative samples that are incorrectly classified as positives. 

A common way to reduce the ROC curve to a single quantity is to calculate the area beneath it \citep{Fawcett2006}. The resulting statistics is called the \emph{area under the curve (AUC)}, and is equivalent to the probability that a classifier will correctly rank a randomly chosen pair of samples $(x, y)$, where $x$ is a sample from $\mathcal{X}$ and $y$ is a sample from $\mathcal{Y}$ \citep{Hanley1982}. The AUC is also equivalent \citep[see][]{Hanley1982,Mason2002} to the popular Mann-Whitney U \citep{Mann1947} and Wilcoxon rank-sum \citep{Wilcoxon1945} statistics, which provide a non-parametric test for differences in the central tendencies of two unpaired samples. It is therefore an appropriate alternative to the two-sample t-test discussed in Section~\ref{sec:unpairedttest}, if the data follow non-Gaussian distributions.

Assuming, without loss of generality, that higher values are indicative for class $\mathcal{X}$, the AUC is given as
\begin{align}
A = Pr(x > y \, | \, x \in \mathcal{X}, y \in \mathcal{Y})\;.
\end{align}
Perfect class separability is denoted by $A=0$ and $A=1$, while chance-level class separability is attained at $A = 0.5$. Thus, a common null hypothesis is $H_0: A = 0.5$.

Assume we have $N_{\mathcal X}$ samples from condition $\mathcal X$ and $N_{\mathcal Y}$ samples from condition $\mathcal Y$. To compute the test statistics, all observations from both conditions are pooled and ranked, beginning with rank one for the smallest value. Defining by $\text{rank}(x_n)$ the rank of $x_n$ (the $n$-th sample from condition $\mathcal{X}$), the Wilcoxon rank-sum statistic for class $\mathcal{X}$ 
is defined as
\begin{align}
W = \sum_{n=1}^{N_{\mathcal X}} \text{rank}(x_n) \;,
\end{align}
while the Mann-Whitney U statistic is given by 
\begin{align}\label{eq:UW_transform}
U = W- \frac{N_{\mathcal X}(N_{\mathcal X} + 1)}{2} \;.
\end{align}
%
Finally, the AUC statistic is given by
\begin{align}\label{eq:AU_transform}
\hat{A} = \frac{U}{N_{\mathcal X} N_{\mathcal Y}} \;.
\end{align}
The exact distributions of $W$, $U$ and $\hat{A}$ under the null hypothesis can be derived from combinatorial considerations \citep{Mann1947,Mason2002}, and critical values for rejecting the null hypothesis can be calculated using recursion formulae \citep{shorack1966recursive}. However, these distributions are approximately normal distributed for samples of moderate size ($N_{\mathcal X} + N_{\mathcal Y} \geq 20$). The mean and variance of Mann-Whitney's $U$ is given by
\begin{align}\label{eq:StdU_H0}
\E_{H_0} (U) &= \frac{N_{\mathcal X} N_{\mathcal Y}}{2} &
{\Var}_{H_0} (U ) &= \frac{N_{\mathcal X} N_{\mathcal Y} (N_{\mathcal X} + N_{\mathcal Y}+1)}{12} \;,
\end{align}
where $\E_{H_0}(\cdot)$ and ${\Var}_{H_0}(\cdot)$ denote expected value and variance under the null hypothesis \citep{Mason2002}.
From Eq.~\eqref{eq:AU_transform}, the mean and variance of the AUC statistic follow as 
\begin{align}\label{eq:StdWA_H0}
\E_{H_0} (\hat{A}) &= \frac{1}{2} & {\Var}_{H_0} (\hat{A} ) &= \frac{ {\Var}_{H_0} (U ) } { N_{\mathcal X}^2 N_{\mathcal Y}^2}  \;.
\end{align}
 Note that this null distribution does not depend on the distribution of the data, and is only based on the assumptions of i.i.d. samples, equal variances of both classes, and that observations are ordinal (that is, it is possible to rank any two observations).

If the null hypothesis is violated (e.g., $A \neq 0.5$), the variances of $U$, $W$, and $\hat{A}$ become data-dependent. The variance for general $A$ can be approximated as \citep{Hanley1982,Greiner2000}
\begin{align}\label{eq:StdAUC_general}
 \widehat{\Var} (\hat{A} ) = \frac{\hat{A} (1 - \hat{A}) + (N_{\mathcal X} - 1) (Q_1 - \hat{A}^2) + (N_{\mathcal Y} - 1) (Q_2 - \hat{A}^2)}{N_{\mathcal X} N_{\mathcal Y}} \;,
\end{align}
where $Q_1 = \hat{A}/(2-\hat{A})$ and $Q_2 = (2 \hat{A}^2)/(1+\hat{A})$. The variances of $U$ and $W$ follow accordingly. 
A statistical test for the null hypothesis can be devised using that %
\begin{align}
z = \frac{\hat{A} - 0.5}{\surd \widehat{\Var}(\hat{A})}
\end{align}
is approximately standard normal distributed for large sample sizes (analogous for $U$ and $W$). 




\subsubsection{Pearson correlation coefficient} \label{sec:correlation}
The \emph{Pearson product-moment correlation coefficient} $\hat{\rho}$ is used when one is interested in the linear dependence of a pair of random variables $(X, Y)$. 
Suppose that for each subject, we have $N$ i.i.d. pairs of observations $(x_{1}, y_{1}), \ldots, (x_{N}, y_{N}) \in \R^2$ with sample mean $(\bar{x}, \bar{y})$. In a neuroimaging context, these pairs could reflect neural activity in two different anatomical structures, or concurrently-acquired neural activity and  behavioral (e.g. response time relative to a stimulus) data. The sample Pearson product-moment correlation coefficient is given by 
\begin{align}
\hat{\rho} = \frac{\sum_{n=1}^{N} (x_{n} - \bar{x}) (y_{n} - \bar{y})}{\sqrt{\sum_{n=1}^{N} (x_{n} - \bar{x})^2}\sqrt{\sum_{n=1}^{N} (y_{j} - \bar{y})^2 }} \;,
\label{eq:pearsoncorr}
\end{align}
where $\hat{\rho} = 1$ denotes perfect correlation, and $\hat{\rho} = -1$ denotes perfect anti-correlation. 
%
%
The null hypothesis of no correlation is given by $H_0: \rho = 0$. Assessing the statistical significance of Pearson correlations can be done using the \emph{Fisher z-transformation} \citep{Fisher1915}, defined as 
\begin{align}\label{eq:FisherTransform}
\zeta(\hat{\rho}) := \frac{1}{2} \ln\left(\frac{1+{\hat{\rho}}}{ 1-{\hat{\rho}}}\right) = \operatorname{arctanh}({\hat{\rho}}) \;.
\end{align}
If $(X, Y)$ has a bivariate normal distribution, then $\zeta(\hat{\rho})$ is approximately normal distributed with mean $ \operatorname{arctanh}(\rho)$ and variance
\begin{align}\label{eq:FisherVariance}
\Var \left( \zeta(\hat{\rho}) \right)= \frac{1}{N - 3} \;.
\end{align}
Therefore the test statistic 
\begin{align}
z = \frac{\zeta(\hat{\rho})}{ \surd \Var \left( \zeta(\hat{\rho}) \right) }
\end{align}
is approximately standard normal distributed. 

The Fisher-transformation is also used when averaging correlations, where the standard approach is to Fisher-transform each individual correlation before computing the average. The reason behind this step is that the distribution of the sample correlation is skewed, whereas the Fisher-transformed sample correlation is approximately normal distributed and thus symmetric \citep[cf.,][]{Silver1987}. Results can be transformed back into a valid Pearson correlation using the inverse transformation
\begin{align}\label{eq:FisherBackTransform}
\hat{\rho} := \frac{{e^{2\zeta(\hat{\rho})} -1}}{e^{2\zeta(\hat{\rho})} +1} = \operatorname{tanh}(\zeta(\hat{\rho})) \; .
\end{align}
The same back transformation can be applied to map confidence intervals derived for $\zeta(\hat{\rho})$ into the Pearson correlation domain. 

Pearson correlation can also be used to derive the \emph{coefficient of determination}, which indicates the proportion of the variance in the dependent variable that is predictable from the independent variable in a linear regression. If an intercept term is included in the regression, the coefficient of determination is given as the square of the Pearson product-moment correlation between the two variables. 
Another strongly related quantity is the point-biserial correlation coefficient, which is used when one variable is dichotomous, i.e., indicates membership in one of two experimental conditions. Pearson correlation is mathematically equivalent to point-biserial correlation if one assigns two distinct numerical values to the dichotomous variable. {Note that Pearson correlation can be seriously biased by outliers. We refer the interested reader to \cite{Pernet2013} for possible remedies.}

\subsubsection{Linear regression coefficients}\label{sec:linear_model}

A multiple linear regression model has the form 
\begin{align}\label{eq:multlin}
y_n = \beta_0 + x_{n,1} \beta_1 + \hdots + x_{n,K} \beta_K + \eta_n \;,
\end{align}
where the dependent variable $y_n, n \in \{1, \hdots, N\}$ is the $n$-th sample, $x_{n,k}, k \in \{1, \hdots, K\}$ are independent variables (or, factors), $\beta_1, \hdots, \beta_K$ are corresponding regression coefficients, $\beta_0$ is an intercept parameter, and $\eta_n$ is zero-mean, uncorrelated noise. In a neuroimaging context, the samples $y_n$ could represent a neural feature such as the activity of a particular brain location measured at various times $n$, while the $x_{n,k}$ could represent multiple factors thought to collectively explain the variability of $y_n$ such as the type of experimental stimulus or behavioral variables. In some fields, such a model is called a {neural encoding} model. It is also conceivable to have the reverse situation, in which the $x_{n,k}$ represent multiple neural features, while the dependent variable $y_n$ is of non-neural origin. This situation would be called neural decoding.

The independent variables $x_{n,k}$ could be either categorial (i.e., multiple binary variables coding for different experimental factors) or continuous. The specific case in which all independent variables are categorial is called analysis of variance (ANOVA). Linear models therefore generalize a relatively broad class of effect size measures including differences between class-conditional means and linear correlations \citep{Poline2012}. 

The most common way to estimate the regression coefficients $\beta_k, k \in \{0, \hdots, K\}$ is ordinary least-squares (OLS) regression. The resulting estimate is also the maximum-likelihood estimate under the assumption of Gaussian-distributed noise. 
Using the vector/matrix notations $\vec{y} = (y_1, \hdots, y_{N})^\top$, $\boldsymbol{\beta} = (\beta_0, \hdots, \beta_K)^\top$, $\boldsymbol{\eta} = (\eta_1, \hdots, \eta_{N})^\top$, $\vec{x}_n = (1, x_{n,1}, \dots, x_{n,K})^\top$, and $\vec{X} = [\vec{x}_1, \hdots, \vec{x}_{N}]^\top \in \R^{N \times (K+1)}$,
Eq.~\eqref{eq:multlin} can be rewritten as $\vec{y} = X \boldsymbol{\beta} + \boldsymbol{\eta}$. The OLS estimate is then given by 
\begin{align}\label{eq:beta}
\hat{\boldsymbol{\beta}} = (\vec{X}^\top\vec{X})^{-1}\vec{X}^\top \vec{y}\;.
\end{align}

The estimated coefficients $\hat{\beta}_k$ can be treated as effect sizes measuring how much of measured data is explained by the individual factors $x_{n,k}$. The null hypothesis for factor $k$ having no explanatory power is $H_0: \beta_k = 0$. The estimated variance of $\hat{\beta}_k$ is 
\begin{align}\label{eq:beta_variance}
\widehat{\Var}(\hat{\beta}_k) = C_{k,k} \;,
\end{align}
where $C = \hat{\sigma}^2_{\!\eta} (\vec{X}^\top\vec{X})^{-1}$ and $\hat{\sigma}^2_{\!\eta} = \frac{1}{N-(K+1)} \sum_{n=1}^{N} (y_n - \hat{\boldsymbol{\beta}}^\top \vec{x}_n)^2 $ is an unbiased estimator of the noise variance. A statistical test for the null hypothesis can be devised using that 
\begin{align}
t = \frac{\hat{\beta}_k} {\surd \widehat{\Var} (\hat{\beta}_k) }
\end{align}
is t-distributed with $N-(K+1)$ degrees of freedom. A similar procedure can be devised for regularized variants such as Ridge regression \citep{hoerl1970ridge}.

\subsection{Nested statistical inference}\label{sec:group_inference}
In the following, our goal is to combine the data of several subjects 
to estimate a population effect and to assess its statistical significance. 
We denote the number of subjects with $S$. The observed effect sizes of each individual subject are denoted by $\hat{\theta}_s$ ($s = 1, \ldots, S$). Other quantities related to subject $s$ are also indexed by the subscript $s$, while the same quantities without subject index denote corresponding group-level statistics.

Two different models may be formulated for the overall population effect. 

\paragraph{Fixed-effect (FE) model}
In the fixed-effect (FE) model, we assume that there is one (fixed) effect size $\theta$ that underlies each subject, that is 
\begin{align}
\theta_1 = \theta_2 = \ldots = \theta_S =: \theta \; .
\end{align}
The observed effect $\hat{\theta}_s$ can therefore be modeled as 
\begin{align}\label{eq:FEModel}
&& \hat{\theta}_s &= \theta + \epsilon_s, & \text{with } \Var(\epsilon_s) &= \sigma_{\!s}^2 \;,  
\end{align}
where $\epsilon_s$ denotes the deviation of the subject's observed effect from the true effect $\theta_s = \theta$. We assume that the noise terms $\epsilon_s$ are independent, zero-mean random variables with subject-specific variance $\sigma_{\!s}^2$.

The null hypothesis tested by a fixed-effect model is that no effect is present in any of the subjects. Thus, $H_0: \theta = \theta_1 = \ldots = \theta_S = \theta_0$, where $\theta_0$ denotes the zero effect. 

\paragraph{Random-effects (RE) model}
In the random-effects (RE) model, the true effect sizes are allowed to vary over subjects. They are assumed to follow a common distribution of effects with mean $\theta$. The observed effect $\hat{\theta}_s$ is modeled as
\begin{align} \label{eq:REModel1}
 && \theta_s &= \theta + \xi_s   &\text{with } \Var(\xi_s) &= \sigma_{\!\text{rand}}^2  \\
 && \hat{\theta}_s &= \theta_s + \epsilon_s &\text{with } \Var(\epsilon_s) &= \sigma_{\!s}^2 \;,  \label{eq:REModel2}
\end{align}
where $\epsilon_s$ denotes the deviation of the subject's observed effect from the true subject-specific effect $\theta_s$, and where $\xi_s$ denotes the deviation of the true subject-specific effect $\theta_s$ from the population effect $\theta$. $\xi_s$ and $\epsilon_s$ are assumed to be zero-mean, independent quantities. The subject-specific variance of $\epsilon_s$ is $\sigma_{\!s}^2$, while the variance of $\xi_s$ is $\sigma_{\!\text{rand}}^2$. For $\sigma_{\!\text{rand}}^2 = 0$, we recover the fixed-effect model. 

The null hypothesis being tested is that the population effect is zero ($H_0: \theta = \theta_0$), while each individual subject-specific effect $\theta_s$ may still be non-zero.\\*

\noindent Besides testing different null hypotheses, fixed-effect and random-effects models assume different variances of the observed effect sizes. In the fixed-effect model, all observed variability is assumed to be within-subject variability 
\begin{align}
\Var(\hat{\theta}_s) = \sigma_{\!s}^2 \;.
\end{align}
The random-effects model additionally accounts for variability between subjects
\begin{align}\label{eq:VarInRE}
\Var(\hat{\theta}_s) = \sigma_{\!s}^2 + \sigma_{\!\text{rand}}^2 \;.
\end{align}
If the data follow a random-effects model, neglecting $\sigma_{\!\text{rand}}^2$ in a fixed-effect analysis leads to an underestimation of the variance. This has negative consequences if we attempt to make inference on the mean population effect ($H_0: \theta = \theta_0$) relying only on a fixed-effect analysis: We may arrive at spurious results, as the underestimated variance leads to p-values that are too low \citep{Hunter2000,Field2003,Schmidt2009}. On the other hand, there is little disadvantage of using a random-effects analysis, even when the data follows a fixed-effect model. As the assumption of a fixed population effect is unrealistic in most practical cases, it is often recommended to carry out random-effects analysis per default \citep{Field2003,Penny2006,Card2011,Monti2011}. 
 
\subsection{Data pooling} \label{sec:pooling}

The most naive approach to conduct group-level inference would be to pool the samples of all subjects, and thus to disregard the nested structure of the data. In electroencephalography (EEG) studies, this approach is sometimes pursued when computing `grand-average' (group-level) waveforms of event-related potentials (ERP) that are elicited by the brain in response to a stimulus. 

Pooling the samples of all subjects may violate the assumption of identically distributed data underlying many statistical tests. Depending on the type of analysis, this may result in an over- or underestimation of the effect size, an over- or underestimation of the effect variance, and ultimately in {over- or underestimated p-values.} 

\begin{figure}
\centering
\sf A \includegraphics[width =\linewidth]{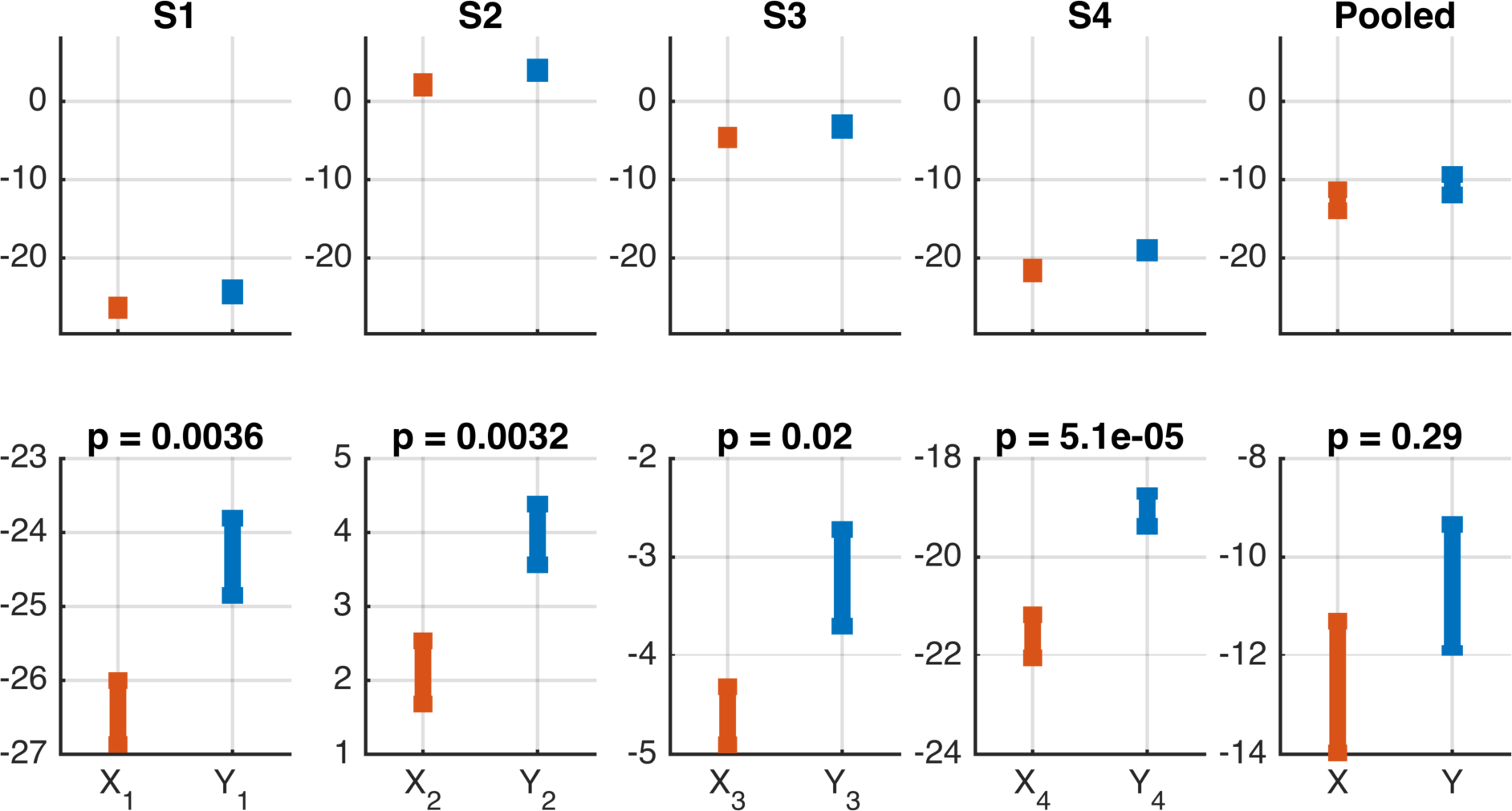} \medskip 

\sf B \includegraphics[width =\linewidth]{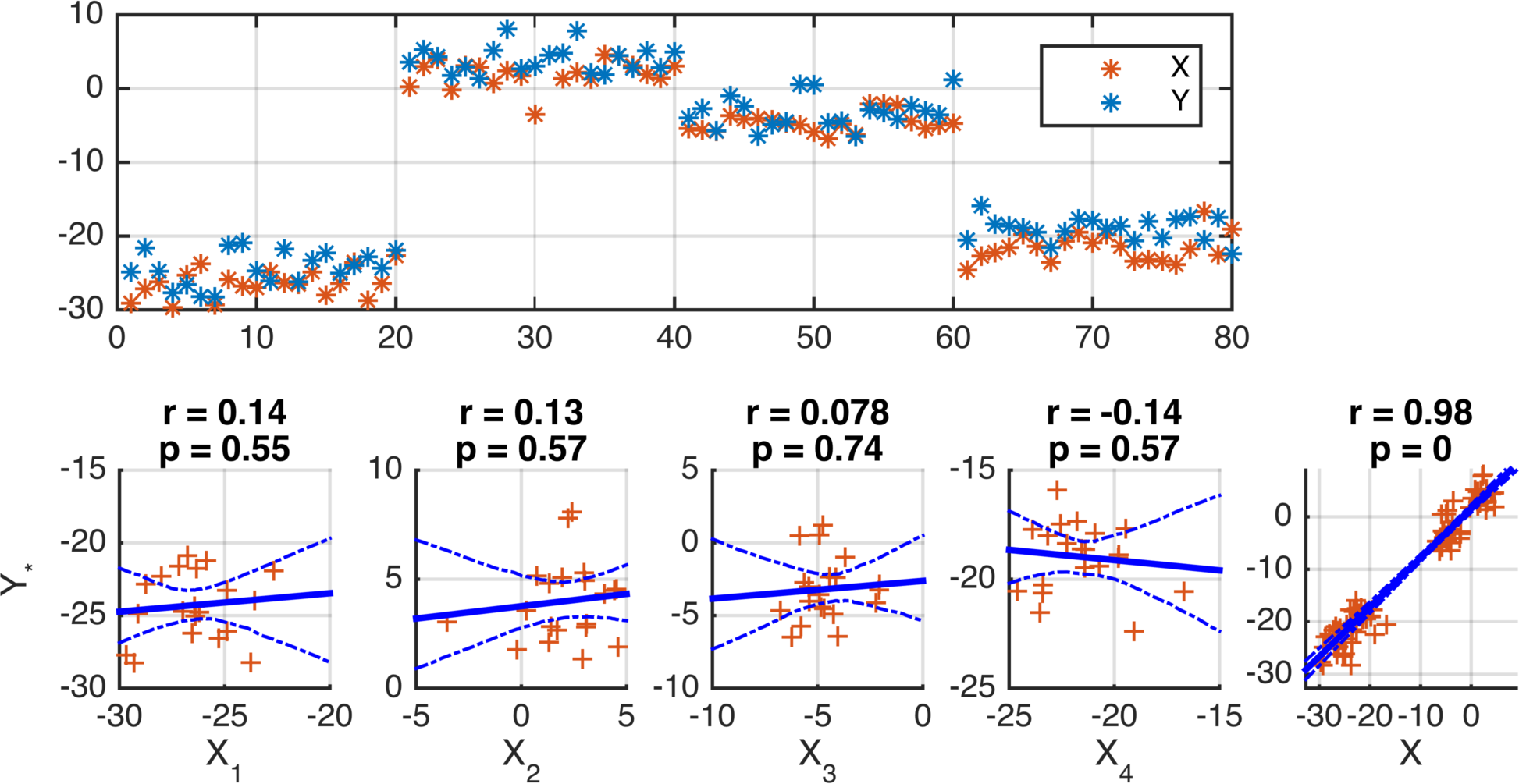} 
\caption{{Wrong conclusions} 
when pooling data with nested structure for statistical testing. Samples were independently drawn for four subjects, $s = 1, \hdots, 4$, and two variables, $X$ and $Y$, according to $x_{n,s} \sim \mathcal{N}(\mu_{s}-1, 4)$, $y_{n,s} \sim \mathcal{N}(\mu_{s}+1, 4)$, where offsets $\mu_{s}$ were drawn independently for each subject from $\mathcal{N}(0, 15^2)$. A) Depiction of the means and standard errors for each subject. A significant difference between means is correctly identified for each subject, but not for the pooled data of all subjects (see lower panels). This is because of the substantial between-subject variance (see upper panels). B) Depiction of the data as a function of sample number (upper panel) and as scatter plots (lower panels). The common subject-specific offsets of $X$ and $Y$ cause strong significant correlation in the pooled data, which is not present in any individual subject, and may be considered spurious. {95\% confidence intervals of the regression line obtained from 1000 Bootstrap samples are marked by dashed blue curves.}
}
\label{fig:pooling_illu}
\end{figure}

The following two examples illustrate the problem. In both cases, two variables, $X$ and $Y$, are modeled for $S=4$ subjects. $N=20$ samples were independently drawn for each subject and variable from Gaussian distributions according to $x_{n,s} \sim \mathcal{N}(\mu_{s}-1, 4)$, $y_{n,s} \sim \mathcal{N}(\mu_{s}+1, 4)$, $s = 1, \hdots, S$, $n=1, \hdots, N$, where the notation $\mathcal{N}(\mu, \sigma^2)$ denotes a Gaussian distribution with mean $\mu$ and variance $\sigma^2$. The subject-specific offsets $\mu_{s}$ were independently drawn from another Gaussian: $\mu_{s} \sim \mathcal{N}(0, 15^2)$. In a practical example, these means may indicate individual activation baselines, which are usually not of interest. Given the generated sample, a difference in the means of $X$ and $Y$ is correctly identified for each subject by Welch's two-sample t-test ($p \leq 0.02$). Because of the substantial between-subject variance, this difference is, however, not significant in the pooled data of all subjects ($p = 0.29$). See Figure~\ref{fig:pooling_illu}~(A) for a graphical depiction. 

A Pearson correlation analysis of the same data correctly rejects the hypothesis of a linear dependence between $X$ and $Y$ for each subject ($|r| \leq 0.14$, $p \geq 0.55$). However, the presence of subject-specific offsets $\mu_s$ causes a strong correlation of $X$ and $Y$ across the pooled data of all subjects ($r = 0.98$, $p \leq 10^{-16}$, see Figure~\ref{fig:pooling_illu}~(B) for a depiction). In many practical cases, this correlation will not be of interest and must be considered spurious.
 
These examples motivate the use of hierarchical approaches for testing data with nested structure, which we introduce below.


\subsection{Naive summary-statistic approach} 
\label{sec:naive}

The simplest variant of the summary-statistic approach ignores subject-specific variances $\sigma_{\!s}^2$, treating subject-level effect sizes $\hat{\theta}_s$ as group-level observations. In this approach, which is somewhat popular in the neuroimaging literature \citep{Holmes1998,Penny2006}, the null hypothesis $H_0: \theta = \theta_0$ is tested based on the $S$ subject-level effect sizes $\hat{\theta}_1, \hdots, \hat{\theta}_S$, which are considered i.i.d.\,. The variance of the mean effect $\hat{\theta} = \nicefrac{1}{S} \sum_{i=s}^S \theta_s$ is estimated as 
\begin{align}
\label{eq:naiveVar}
\widehat{\Var}(\hat{\theta}) = \frac{1}{S-1}\sum_{s=1}^{S} (\hat{\theta}_s - \hat{\theta})^2 \;,
\end{align}
which is an unbiased estimate of $\Var(\hat{\theta})$ even if variances $\sigma_{\!s}^2$ vary across subjects \citep{Mumford2009}. If the $\theta_s$ are normal distributed (for example, because they represent the means of normal distributed or many quantities), the test statistic
\begin{align}
t =\frac{\hat{\theta} - \theta_0}{\surd \widehat{\Var}(\hat{\theta})}
\end{align}
is t-distributed with $S-1$ degrees of freedom. This is the standard one-sample t-test applied to the individual effect sizes $\theta_1, \ldots, \theta_S$.

The naive summary-statistic approach is valid both under the fixed-effect and random-effects models \citep{Mumford2009}. Its statistical power is, however, limited due to two factors. First, it assigns equal importance to each subject. This is sub-optimal if subject-level variances $\sigma_{\!s}^2$ vary across subjects (for example, because of different amounts of recorded data). In this case, a weighting scheme taking into account subject-level variances is optimal (see Section~\ref{sec:invvar}). Second, the approach does not make use of all the available data, as only the group level data is used to  estimate the variance $\widehat{\Var}(\hat{\theta})$ through Eq.~\eqref{eq:naiveVar}. However, even if subject-level variances  $\sigma_{\!s}^2$ are constant across subjects, it is beneficial to make use of their estimates (see Section~\ref{sec:equalvar}).

Both issues are addressed by the sufficient-summary statistic approach described in the next section. An empirical comparison of the statistical power of both approaches on simulated data is provided in Section~\ref{sec:simulations}.



\subsection{Sufficient-summary-statistic approach}
\label{sec:anylincomb}
If estimates of the variances $\Var(\hat{\theta}_s)$ of the subject-level effect sizes $\hat{\theta}_s, s = \{1, \hdots, S\}$ can be obtained, this gives rise to a more powerful summary-statistic approach compared to the naive approach outlined above. {To this end, we estimate the group-level effect size estimate $\hat{\theta}$ as a convex combination 
\begin{align} \label{eq:weightedaverage}
\hat{\theta} := \frac{\sum_{s=1}^S \alpha_s \hat{\theta}_s}{\sum_{s=1}^S \alpha_s}
\end{align}
of the subject-level effect size estimates $\hat{\theta}_s$ with non-negative weights $\alpha_s, \; s \in \{1, \ldots , S\}$. Under the assumption that the $\hat{\theta}_s$ are unbiased and statistically independent of the weights $\alpha_s$, $\hat{\theta}$ is also unbiased (has expectation $\E(\hat{\theta}) = \theta$), as the denominator of Eq.~\eqref{eq:weightedaverage} ensures that the weights sum to one.}
Importantly, with the exception of the coefficient of determination discussed in Section~\ref{sec:correlation}, all effect size measures discussed in this paper are unbiased estimators of the corresponding population effects. The variance of $\hat{\theta}$ defined in Eq.~\eqref{eq:weightedaverage} is given by
\begin{align}
\Var(\hat{\theta}) = \frac{ \sum_{s=1}^S \alpha_s^2 \Var(\hat{\theta}_s) }{ \left(\sum_{s=1}^S \alpha_s \right)^2 } \;.
\label{eq:weightedvariance}
\end{align}
If the $\hat{\theta}_s$ are normal distributed (for example, because they represent the means of normal distributed  or many quantities), the weighted mean $\hat{\theta}$ is also normal distributed. According to the central limit theorem, this is also approximately the case if the $\hat{\theta}_i$ are not normal distributed but the number of subjects $S$ is large. In both cases, we can test the null hypothesis $H_0 : \theta = \theta_0$ using that the test statistic 
\begin{align}\label{eq:z_weightedmean}
z = \frac{\hat{\theta} - \theta_0}{\surd\Var(\hat{\theta})}
\end{align}
is standard normal distributed.

The variances $\Var(\hat{\theta}_s)$ typically need to be estimated, as the exact population values are unknown. As any estimate $\widehat{\Var}(\hat{\theta})$ integrates information from all samples of all $S$ subjects, it can be considered a fairly accurate estimate, justifying the use of a z-test even when we replace $\Var(\hat{\theta})$ by it's estimate $\widehat{\Var}(\hat{\theta})$ in Eq.~\eqref{eq:z_weightedmean} \citep{Borenstein2009,Card2011}. Sometimes, however, the more conservative t-distribution with $S-1$ degrees of freedom is assumed for $z$ \citep{Thirion2007,Jackson2010}.

\subsubsection{Equal weighting}
\label{sec:equalvar}
The z-test introduced in Eq.~\eqref{eq:z_weightedmean} is valid regardless of the choice of the non-negative weights $\alpha_s, \; s \in \{1, \ldots , S\}$ {as long as these weights are statistically independent of the corresponding effect size estimates.} One popular choice is to assign equal weights 
\begin{align}\label{eq:equalweights}
\alpha_1 = \hdots = \alpha_S = \frac{1}{S}
\end{align}
to all subjects, such that $\hat{\theta}$ becomes the arithmetic mean of the $\hat{\theta}_s$. 
This procedure is similar to the naive summary-statistic approach introduced in Section~\ref{sec:naive} in that both approaches assign equal importance to each subject-level effect size. However, it differs in the way the variance is estimated, and in terms of the distribution that is assumed for the test statistic.  For the naive summary-statistic approach, variances are estimated through Eq.~\eqref{eq:naiveVar} using the $S$ data points on the group-level only. The equal-weighting approach instead uses the subject-level variances. That is, following Eq.~\eqref{eq:weightedvariance}: 
\begin{align}\label{eq:equalweightsvariance}
\widehat{\Var}(\hat{\theta}) & = 1/S^2 \sum_{s=1}^S \widehat{\Var}(\hat{\theta}_s) \;.
\end{align} 
If the individual $\widehat{\Var}(\hat{\theta}_s)$ are unbiased, both methods yield an unbiased estimate of the variance $\Var(\hat{\theta})$. But the variance of this variance estimate is typically smaller for the equal variance weighting approach, because it makes use of all the available data. This more accurate estimate means that the test statistic is approximately normal distributed rather than t-distributed with $S-1$ degree of freedoms. This translates into a power gain, as illustrated in the simulation presented in Section~\ref{sec:simulations}. However, estimating the between-subject variance for a random-effects model is not straightforward, and also may introduce biases and variability (see Section~\ref{sec:sigmarand}).

\subsubsection{Inverse-variance weighting}
\label{sec:invvar}
Interestingly, the choice of equal weights is suboptimal in terms of obtaining a group-level effect size estimate $\hat{\theta}$ with minimal variance. It is generally desirable to minimize the variance of the weighted average, as unbiased estimators with smaller variance achieve a lower mean squared error (MSE), and lead to more powerful statistical tests. The minimum-variance estimate is obtained by weighting each subject-level effect size proportional to the inverse of its variance using weights 
\begin{align}\label{eq:invvarweights}
\alpha_s = \frac{1 }{ \Var(\hat{\theta}_s) }\;.
\end{align}
This result is consistent with the intuition that less precise $\theta_s$ should have a lower impact on the overall estimate than those that are estimated with high confidence. Inserting into Eq.~\eqref{eq:weightedvariance}, we obtain the optimal value
\begin{align}\label{eq:invvarweightsvariance}
\Var(\hat{\theta})= \frac{1}{\sum_{s=1}^S 1 /\Var(\hat{\theta}_s) } = \frac{1}{\sum_{s=1}^S \alpha_s} \;.
\end{align}
{Note, however, that, by using data-dependent weights, the inverse-variance-weighting approach may not always result in unbiased group-level effect size estimates. The potential implications of correlations between individual subject-level effect sizes and their variances are demonstrated in Section~\ref{sec:bias_example} and further discussed in Section~\ref{sec:invvar_limitation}.}


\subsubsection{Estimation of between-subject variance}
\label{sec:sigmarand}

To perform inverse-variance weighting under the random-effects model, the between-subjects variance $\sigma_{\!\text{rand}}^2$ needs to be estimated in order to obtain the total subject-wise variance $\Var(\hat{\theta}_s) = \sigma_{\!s}^2 + \sigma_{\!\text{rand}}^2$. Several iterative and non-iterative alternative methods have been proposed \citep{Worsley2002,Guolo2015,Veroniki2016}. A popular and easy-to-implement approach is the non-iterative procedure proposed by \cite{Dersimonian1986}. For a given estimate $\hat{\sigma}_{\!s}^2$ of the within-subject variances (which can be obtained using the procedures discussed in Section~\ref{sec:effectsizes}), and for fixed-effect quantities
\begin{align}\label{eq:DersimonianLaird0}
\alpha_s^{FE} = \frac{1}{\hat{\sigma}_{\!s}^2}, \qquad \hat{\theta}^{FE} = \frac{\sum_{s=1}^S \alpha_s^{FE} \hat{\theta}_s}{\sum_{s=1}^S \alpha_s^{FE} } \;,
\end{align}
the between-subject variance $\sigma^2_{\!\text{rand}}$ according to \cite{Dersimonian1986} is estimated as 
\begin{align}\label{eq:DersimonianLaird}
\hat{\sigma}^2_{\!\text{rand}} = \max \left\{ 0, \frac{\sum_{s=1}^S \alpha_s^{FE} (\hat{\theta}_s - \hat{\theta}^{FE})^2 - S + 1}{\sum_{s=1}^S \alpha_s^{FE} - \sum_{s=1}^S (\alpha_s^{FE})^2 / \sum_{s=1}^S \alpha_s^{FE} } \right\} \;.
\end{align}
{where $\alpha_s^{FE}$ and $\hat{\theta}^{FE}$ are the fixed effect quantities defined in Eq.~\eqref{eq:DersimonianLaird0}.}

As this estimate may be quite variable for small sample sizes, the resulting p-values may become too small when the number of subjects $S$ is small \citep{Brockwell2001,Guolo2015}. On the other hand, the truncation of the estimated variance to zero introduces a positive bias; that is, $\sigma^2_{\!\text{rand}}$ (and thus, p-values) are generally over-estimated \citep{Rukhin2013}. In summary, the Dersimonian and Laird approach is acceptable for a moderate to large number of subjects \citep{Jackson2010,Guolo2015}, and is the default approach in many software routines in the meta-analysis community \citep{Veroniki2016}. 

After $\hat{\sigma}^2_{\!\text{rand}}$ has been calculated, the random-effects quantities are finally computed as
\begin{align}
\alpha_s^{RE} = \frac{1}{\hat{\sigma}_{\!s}^2 + \hat{\sigma}_{\!\text{rand}}^2}, \qquad \hat{\theta}^{RE} = \frac{\sum_{s=1}^S \alpha_s^{RE} \hat{\theta}_s}{\sum_{s=1}^S \alpha_s^{RE} } \; .
\end{align}

%


\subsubsection{Algorithm}

The sufficient-summary-statistic approach is summarized in Algorithm~\ref{algo:invVar}. 
First, the subject-level effect sizes $\theta_s$ and the within-subject variances ${\sigma}_{\!s}^2, s \in \{1, \hdots, S\},$ are estimated based on the available subject-wise data samples. 
Second, if random-effects are assumed, the correlation between $\theta_s$ and ${\sigma}_{\!s}, s \in \{1, \hdots, S\}$ across subjects is assessed, preferably using a robust measure such as Sperman's rank correlation. 
Third, the between-subject variance ${\sigma}_{\!\text{rand}}^2$ is estimated as outlined in Section~\ref{sec:sigmarand} (unless a fixed-effect model can reasonably be assumed). The variance of a subject's estimated effect $\hat{\theta}_s$ around the population effect $\theta$ is calculated as the sum of the within-subject measurement error variance $\sigma_{\!s}^2$ and the between-subject variance $\sigma_{\!\text{rand}}^2$ (cf. Eq.~\eqref{eq:VarInRE}). 
Fourth, the estimated population effect $\hat{\theta}$ is calculated as the weighted average of the subjects effects. If a fixed-effect is assumed, or if no correlation between effect size and variances estimates has been found in the random effects setting, weights $\alpha_s, s = 1, \hdots, S$ are set to the inverse of the estimated subject-level variances as outlined in Section~\ref{sec:invvar}. If a correlation between subject-level effect sizes and standard deviations has been found, it is instead advisable to use equal weights for all subjects (Section~\ref{sec:equalvar}) or weights that are proportional to the subjects' sample sizes \citep{marin2010weighting}. Given the weights $\alpha_s$, the variance of the variance of the population effect can be calculated either using the general formula given by Eq.~\eqref{eq:weightedvariance} or specific versions derived for equal and inverse-variance weighting schemes in Eqs.~\eqref{eq:equalweightsvariance} and \eqref{eq:invvarweightsvariance}. Finally, the estimated mean effect is subjected to a z-test as introduced in Eq.~\eqref{eq:z_weightedmean}.

Different effect sizes and their corresponding variances have been discussed in Section~\ref{sec:effectsizes}. 
With the exception of the Pearson correlation coefficient, these measures can be directly subjected to the inverse-variance-weighting approach. That is, $\hat{\theta}_s$ and $\hat{\sigma}_{\!s}^2$ for the mean difference are given in Eqs.~\eqref{eq:meanDifference} and \eqref{eq:meanDifferenceVariance}, for the AUC in Eqs.~\eqref{eq:AU_transform} and \eqref{eq:StdAUC_general}, and for linear regression coefficients in Eqs.~\eqref{eq:beta} and \eqref{eq:beta_variance}. As discussed in Section~\ref{sec:correlation}, it is, however, beneficial to transform correlation coefficients $\hat{\rho}_s$ into approximately normal distributed quantities with known variance prior to averaging across subjects. We can proceed with the application of the sufficient-summary-statistic approach just as outlined before, treating the transforms $\zeta(\hat{\rho}_s)$ given in Eq.~\eqref{eq:FisherTransform} rather than the $\hat{\rho}_s$ as effect sizes. The resulting population effect can be transformed back into a valid Pearson correlation using the inverse transformation described in Eq.~\eqref{eq:FisherBackTransform}. 

\begin{algorithm}[htbp]
\begin{algorithmic}
\STATE {\bf Step 1: Within-subject analysis}
\bindent
 \FORALL {Subjects $s = 1 \ldots S$}
  \STATE Estimate effect size $\hat{\theta}_s$ and its variance $\hat{\sigma}_{\!s}^2$
\ENDFOR 
\eindent
\STATE {\bf Step 2 : Correlation between effect size and variance}
\bindent
\STATE Random effects setting: test $H_0^{\text{corr}}: \rho_{\theta_s, \sigma_{\!s}} = 0$
\STATE Fixed effect setting: accept $H_0^{\text{corr}}$
\eindent

\STATE {\bf Step 3: Between-subject variance $\sigma_{\!\text{\normalfont rand}}^2$}
\bindent
\STATE Random effects setting: use, e.g., Eq. \eqref{eq:DersimonianLaird0}-\eqref{eq:DersimonianLaird}
\STATE Fixed effect setting: $\hat{\sigma}_{\!\text{rand}}^2 \gets 0$
\eindent
\STATE {\bf Step 4: Population mean effect and variance}
\bindent
\FORALL {Subjects $s = 1 \ldots S$}
\IF {$H_0^{\text{corr}}$ is accepted}
\STATE Perform inverse-variance weighting:
\STATE $\alpha_s \gets 1 / ( \hat{\sigma}_{\!s}^2 + \hat{\sigma}_{\!\text{rand}}^2 )$ 
\STATE $\hat{\theta} \gets \sum_{s=1}^S \alpha_s \hat{\theta}_s / \sum_{s=1}^S \alpha_s $ 
\STATE $\widehat{\Var}(\hat{\theta}) \gets 1 / \sum_{s=1}^S \alpha_s $
\ELSE
\STATE $\alpha_s \gets 1 / S$ (equal weighting) \textbf{\;or\;} $\alpha_s \gets N_s / \sum_{s=1}^S N_s$ (sample-size weighting)
\STATE $\hat{\theta} \gets \sum_{s=1}^S \alpha_s \hat{\theta}_s $ 
\STATE $\widehat{\Var}(\hat{\theta}) \gets \sum_{s=1}^S \alpha_s^2 ( \hat{\sigma}_{\!s}^2  + \hat{\sigma}_{\!\text{rand}}^2)$
\ENDIF
\ENDFOR 
\eindent
\STATE {\bf Step 5: Statistical inference ($H_0: \theta = \theta_0$)}
\bindent
\STATE $z \gets (\hat{\theta} - \theta_0) / \surd\widehat{Var}(\hat{\theta})$ 
\STATE $z$ is approximately standard normal distributed
\STATE $\Rightarrow$ Reject $H_0$ at $0.05$ level if $|z| > 1.96$ 
\eindent
\end{algorithmic}
\caption{Sufficient-summary-statistic approach}
\label{algo:invVar}
\end{algorithm}

\subsection{Stouffer's method of combining p-values}
\label{sec:stouffer}

A general approach for combining the results of multiple statistical tests is Stouffer's method \citep{Stouffler1949,Whitlock2005}. For a set of independent tests of null hypotheses $H_{0, 1}, \hdots, H_{0, S}$, Stouffer's method aims to determine whether all individual null hypotheses are jointly to be accepted or rejected, or, in other words, if the global null hypothesis $H_0: \text{(} \forall s: H_{0, s}$ is true) is true. In general, the individual $H_{0,s}$ may not necessarily refer to the same effect size or even effect size measure, and the p-values for each individual hypothesis may be derived using different test procedures including non-parametric, bootstrap- or permutation-based tests. In the present context of nested data, Stouffer's method can be used to test group-level null hypotheses in the fixed-effect setting, i.e., the absence of an effect in all $S$ subjects of the studied population. 

Denote with $H_{0,s} : \theta_s = \theta_0$ the null hypothesis that there is no effect in subject $s$, and with $p_s$ the \emph{one-tailed} p-value of an appropriate statistical test for $H_{0,s}$. If the null hypothesis is true, $p_s$ is uniformly distributed between $0$ and $1$ \citep[see][for an illustration]{Murdoch2012}. Therefore, the one-tailed p-values $p_s$ can be converted into standard normal distributed z-scores using the transformation
\begin{align}\label{eq:Stoufferz}
z_s := F_{z_{0,1}}^{-1} (p_s) \;,
\end{align}
where $F_{z_{0,1}}^{-1}$ denotes the inverse of the standard normal cumulative distribution function. For Gaussian-distributed subject-level effect sizes with known variance, this step can be carried out more directly using
\begin{align}
{z}_s = \frac{\hat{\theta}_s - \theta_0}{ \surd \Var(\hat{\theta}_s)} \;.
\end{align}
The cumulative test statistic 
\begin{align}\label{eq:Stouffer_grouplevel}
z = \frac{1}{\surd{S}} \sum_{s = 1}^S z_s
\end{align}
follows the standard normal distribution, which can be used to derive a p-value for the group-level $H_0$.

Notice that Stouffer's method as outlined above is applied to \textit{one-tailed} p-values. However, testing for the presence of an effect often requires a two-tailed test. In this case, it is important to take the direction of the effect in different subjects into account. We cannot simply combine two-tailed tests -- a positive effect in one subject and a negative effect in another subject would be seen as evidence for an overall effect, even though they cancel each other out. However, the direction of the effect can be determined post-hoc. To this end, one-tailed p-values for the same direction are calculated for each subject and combined as outlined in Eqs.~\eqref{eq:Stoufferz} and \eqref{eq:Stouffer_grouplevel} into a group-level one-tailed p-value $p_1$. The group-level two-tailed p-value is then obtained as $p_2 = 2 \cdot \min(p_1, 1-p_1)$ (see Eq.~\eqref{eq:one-tailed}-\eqref{eq:two-tailed}) \citep{Whitlock2005} .

\section{Simulations}
\label{sec:simulations}

\begin{figure*}
\centering
\includegraphics[width =0.43\textwidth]{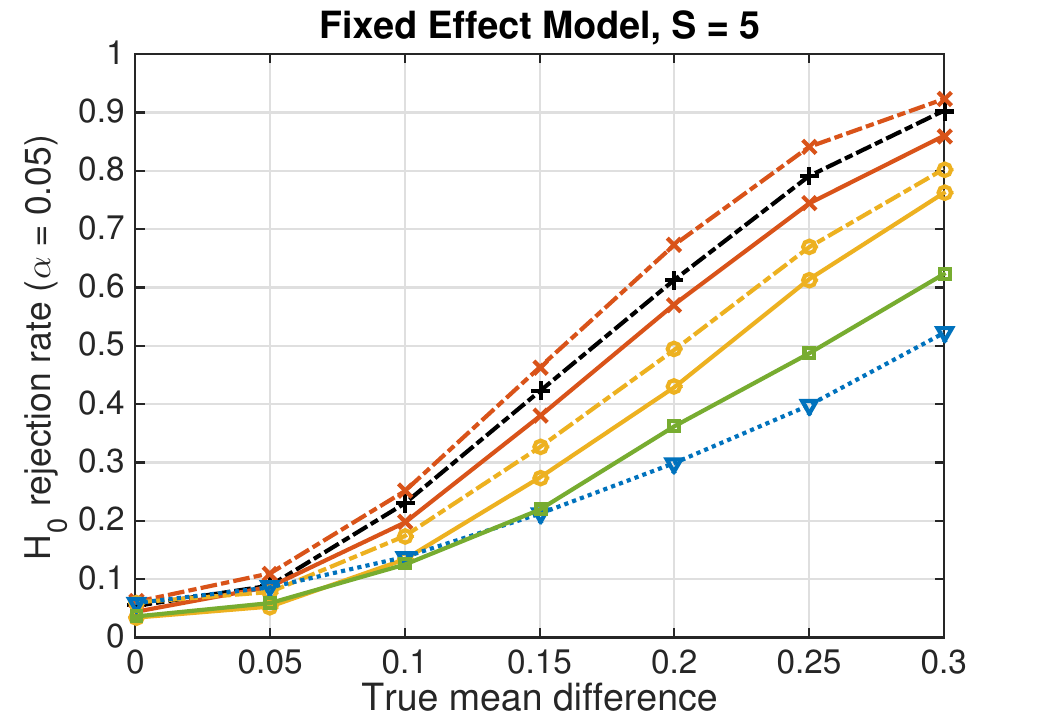} 
\includegraphics[width =0.43\textwidth]{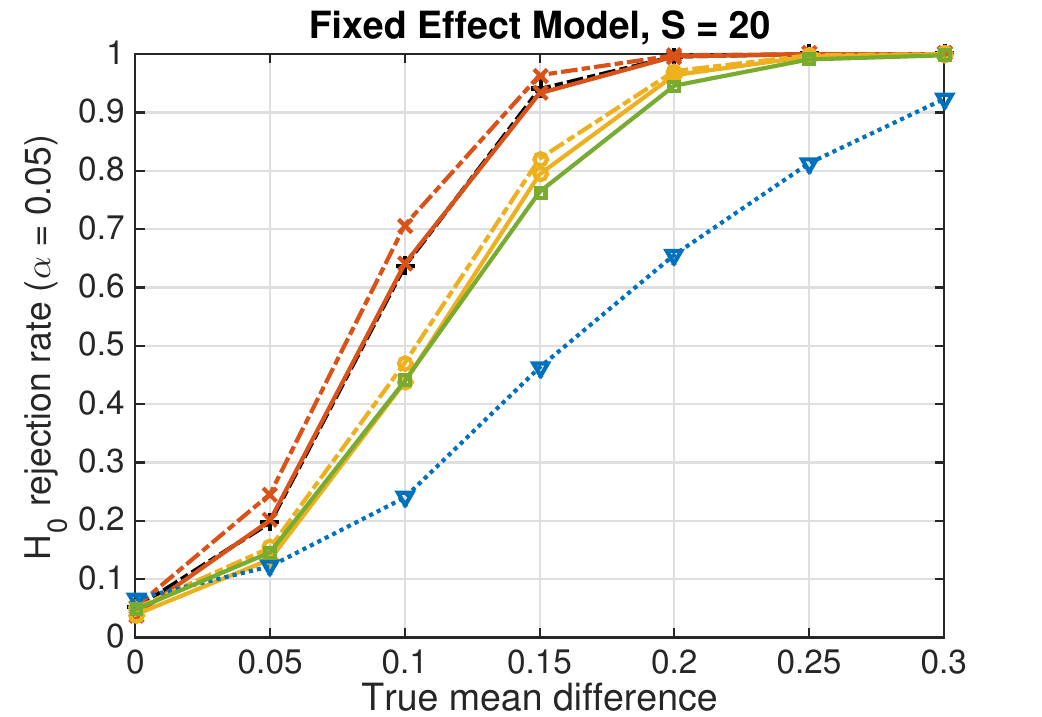}  \\ \vspace{2ex}
\includegraphics[width =0.43\textwidth]{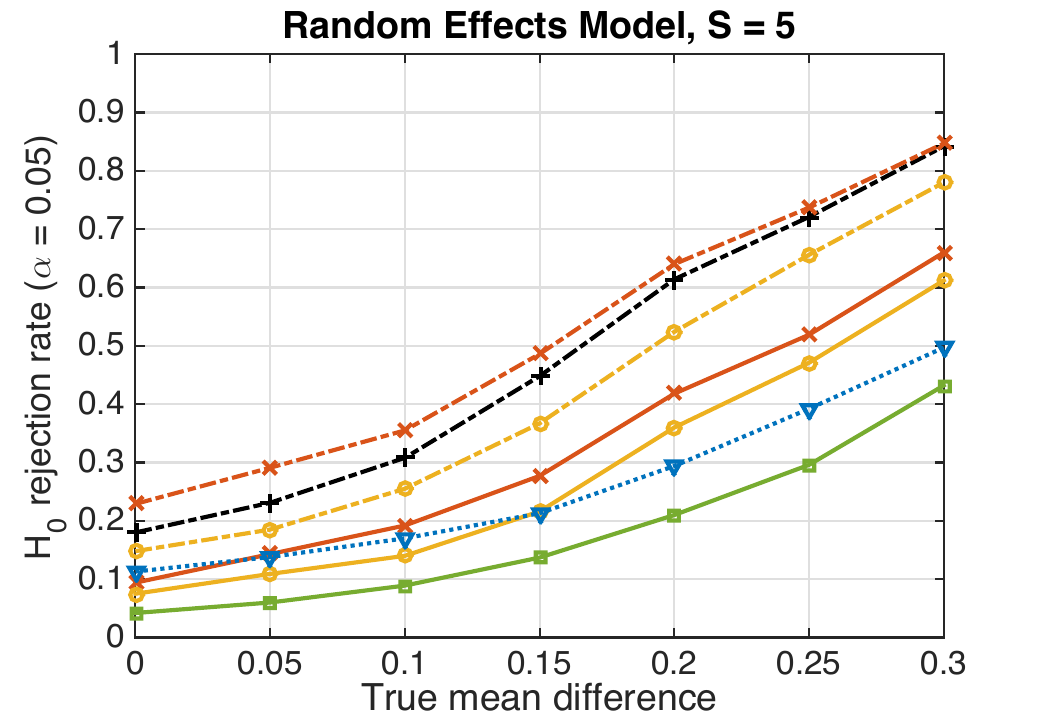} 
\includegraphics[width =0.43\textwidth]{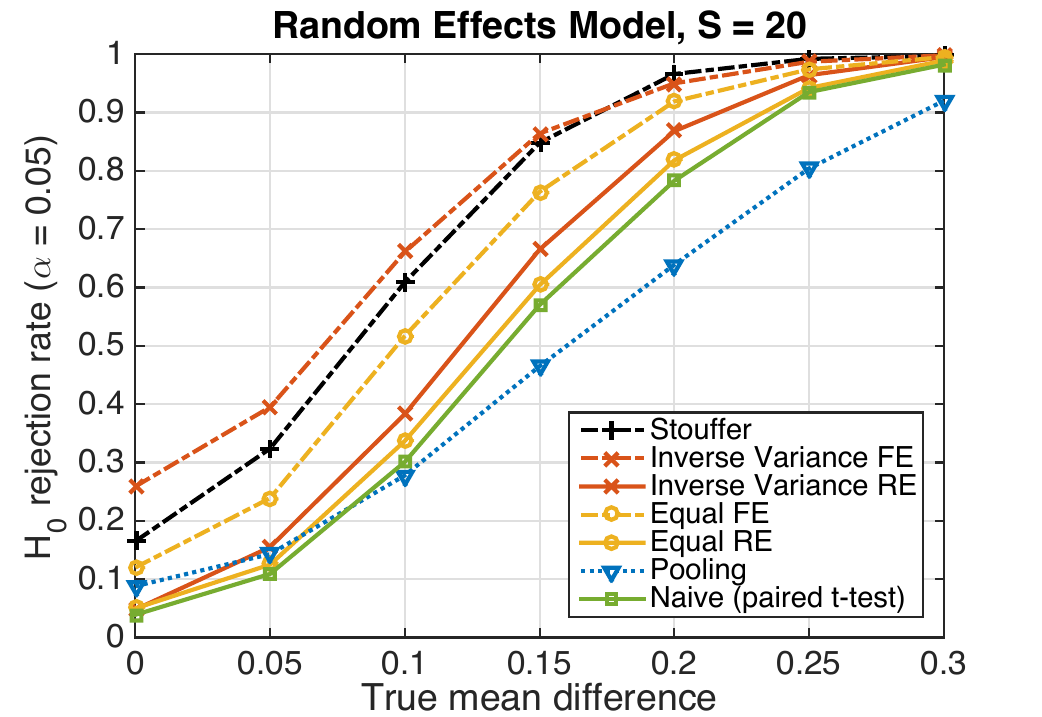} \\ 
\caption{The probability of rejecting the null hypothesis $H_0: d = 0$ as a
function of the true mean difference $d$ of Gaussian-distributed simulated data from $S=5$ resp. $S=20$ subjects. Top: data following a fixed-effect model. Bottom: data following a random-effects model.}
\label{fig:simul1}
\end{figure*}

In the following, we present a set of simulations, in which we compare the statistical approaches reviewed above to test for a difference between two class-conditional means in artificial data. We consider two conditions $\mathcal{X}$ and $\mathcal{Y}$ with true means $\mu_{\mathcal{X}}$ and $\mu_{\mathcal{Y}}$ and class-conditional mean difference $d = \mu_{\mathcal{Y}} - \mu_{\mathcal{X}}$. We want to test the null hypothesis $H_0 : \mu_{\mathcal{X}} = \mu_{\mathcal{Y}}$ or, equivalently, $H_0 : d = 0$. The following scenarios are investigated: 1) The data are generated either within a fixed-effect or a random-effects model. 2) The data are generated from either a Gaussian or a non-Gaussian distribution. In each scenario, we compare the methods' abilities to reject the null hypothesis when we vary the true class-conditional mean difference $d$. 

Data for $S = 5$ or $S=20$ subjects $s, s \in \{1, \ldots, S \},$ were generated as follows. First, subject-specific class-conditional mean differences $d_s$ were sampled according to 
\begin{align*}
d_s = d + \xi_s \;, \qquad \xi_s \sim \mathcal{N}(0, \sigma_{\!\text{rand}}^2)  \;,
\end{align*}
where $\sigma_{\!\text{rand}}^2$ is the between-subject variance. For the fixed-effect model, we set $\sigma_{\!\text{rand}} = 0$, while for the random-effects model, we set $\sigma_{\!\text{rand}} = 0.2$.

We then sampled $N_{s,\mathcal{X}}$ data points for condition $\mathcal{X}$ and $N_{s,\mathcal{Y}}$ data points for condition $\mathcal{Y}$ from Gaussian distributions with variance $v^2_s$ and class-conditional means $\mu_{s,\mathcal{X}}$ and $\mu_{s,\mathcal{Y}} = \mu_{s,\mathcal{X}} + d_s$, respectively. A separate set of samples was drawn from non-Gaussian F(2,5)-distributions adjusted to have the same class-conditional means and variance. The number of data points, $N_{s,\mathcal{X}}$ and $N_{s,\mathcal{Y}}$, the class-conditional means, $\mu_{s,\mathcal{X}}$ and $\mu_{s,\mathcal{Y}}$, and the variance, $v^2_s$, were randomly drawn for each subject such that $v_s$~is uniformly distributed between 0.5 and 2, $N_{s,\mathcal{X}}$ and $N_{s,\mathcal{Y}} \in \mathbb{N}$ are uniformly distributed between 50 and 80, and the true mean of class $\mathcal{X}$, $\mu_{s,\mathcal{X}}$, is uniformly distributed between -3 and 3. In each scenario, the true class-conditional mean difference, $d$, was varied across $\{0, 0.05, 0.1, 0.15, 0.2, 0.25, 0.3\}$. 
 
All experiments were repeated 1000 times with different samples drawn from the same distributions. We report the \textit{$H_0$ rejection rate}, which is the fraction of the test runs in which the null hypothesis was rejected. When the null hypothesis is true ($d = 0$), the $H_0$ rejection rate is identical to the error or false positive rate of the statistical tests under study. In the converse case, in which the null hypothesis is false ($d \neq 0$), the rejection rate determines the power of the test. All statistical tests were performed at significance level $\alpha = 0.05$. An ideal test would thus obtain a $H_0$ rejection rate of $0.05$ when the null hypothesis is true ($d=0$), and a rejection rate of 1 otherwise. The higher the $H_0$ rejection rate in the presence of an effect ($d \neq 0$), the higher is the power of a test. However, if the null hypothesis is true, a $H_0$ rejection rate greater than $\alpha$ indicates the occurrence of spurious findings beyond the acceptable $\alpha$-level.

\begin{figure*}
\centering
\includegraphics[width =0.37 \linewidth]{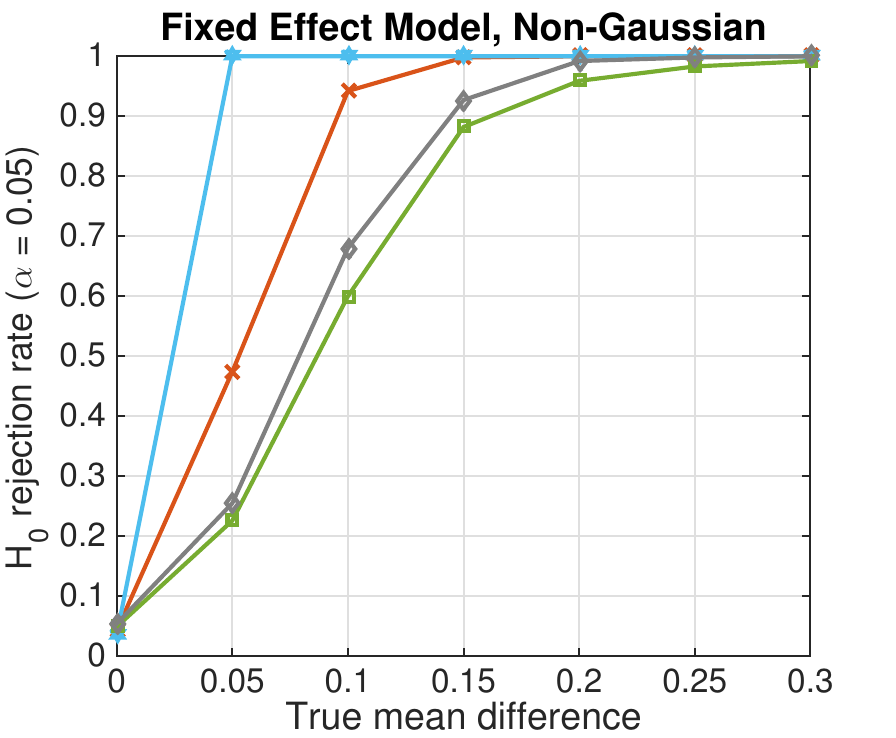} 
\includegraphics[width =0.37 \linewidth]{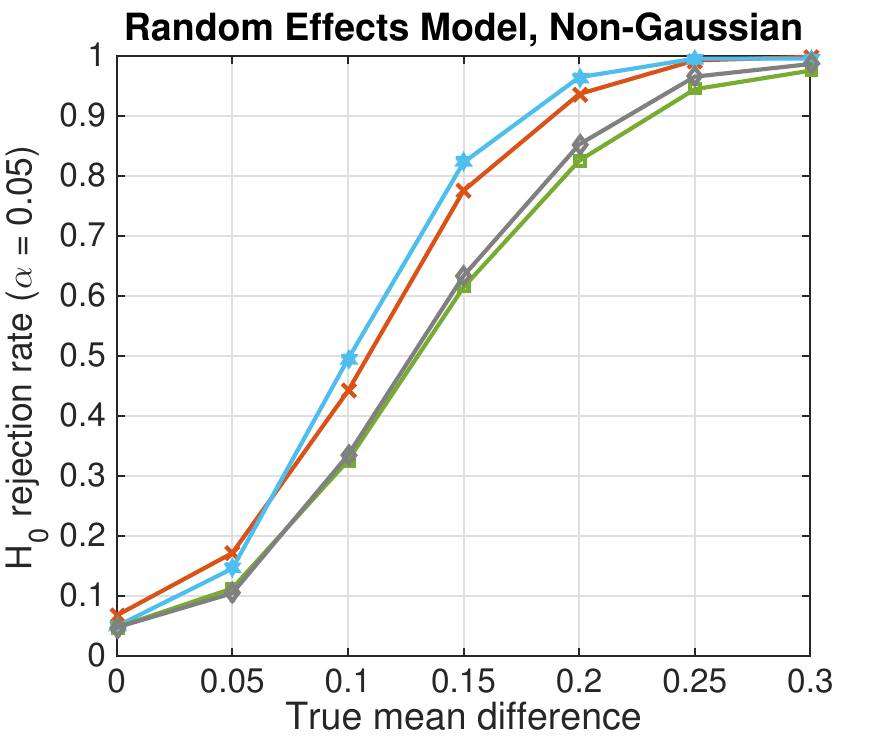} 
\includegraphics[width =0.37 \linewidth]{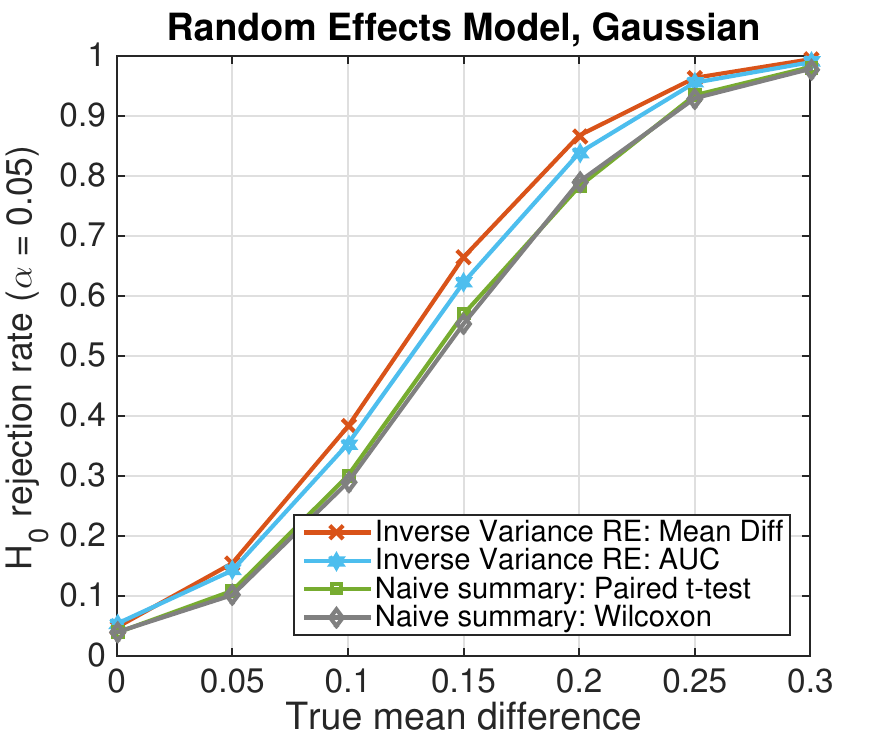} 
\caption{The probability of rejecting the null hypothesis $H_0: d = 0$ as a
function of the true mean difference $d$ in simulated data of $S=20$ subjects. Top Right: Non-Gaussian data from a fixed effect model. Top left: Non-Gaussian data from a random effects model. Bottom: Gaussian data from a random effects model. }
\label{fig:simul2}
\end{figure*}

\subsection{Simulation 1: Fixed effect vs. random effects}

Figure~\ref{fig:simul1} depicts the results achieved by the tested statistical procedures in the fixed-effect (top row) and random-effects (bottom row) scenarios for Gaussian-distributed data, using data from $S=5$ and $S=20$ subjects. The `pooling' approach consists of pooling the samples of all subjects and performing one two-sample t-test (cf. Section~\ref{sec:pooling}). `Naive (paired t-test)' refers the naive summary-statistic approach, in which each subject's mean difference is treated as an observation for a group-level paired t-test (cf. Section~\ref{sec:naive}). Four variants of the sufficient-summary-statistic approach are considered (cf. Section~\ref{sec:anylincomb}). These variants differ in assuming either random effects (RE) or one fixed effect (FE), and in using either the inverse-variance-weighting scheme (Eq.~\eqref{eq:invvarweights}) or equal weights (Eq.~\eqref{eq:equalweights}). `Stouffer' finally refers to using Stouffer's method to combine p-values obtained from subject-level two-sample t-tests  (cf. Section~\ref{sec:stouffer}). Note that all group-level tests are carried out two-tailed.

In line with our previous considerations, data pooling yielded very low power in the presence of an effect both under the fixed-effect and random-effects models. The highest power is achieved in both cases by the inverse-variance-weighted sufficient-summary-statistic approach, followed by Stouffer's method, the sufficient-summary-statistic approach using equal weights, and the paired t-test. 

Considerable differences are observed between the fixed-effect and random-effects settings. For data following the fixed-effect model, the fixed-effect variants of the sufficient-summary-statistic approach display only a negligible advantage over their random-effects counterparts, indicating that the latter succeed in estimating the between-subject variance to be zero. Moreover, in the case of equal class means, all approaches achieve a false positive rate close to the expected value of $\alpha=0.05$. 

The situation is different for data following a random-effects model. Here, the fixed-effect variants of the sufficient-summary-statistic approach as well as Stouffer's method and the pooling approach display false positive rates that are between two and five times higher (26\%) than what would acceptable under the null hypothesis. This problem is substantially alleviated by the random-effect variants of the sufficient-summary-statistic approach. Nevertheless, when data is only available from $S=5$ subjects, the null hypothesis is still rejected too often ($9\%$ for inverse-variance weighting). This is due to the variability in the estimate of the between-subject variance $\sigma_{\!\text{rand}}^2$ (cf. Section~\ref{sec:sigmarand}). When $S=20$ subjects are available, the expected false positive rate of $\alpha = 0.05$ is achieved.

The naive summary-statistic approach (paired t-test of subject-wise means) achieves the expected false positive rate of $0.05$ regardless of the number of subjects, and therefore  represents a valid statistical test also in the random-effects setting.

\subsection{Simulation 2: Gaussian vs. non-Gaussian}

Figure~\ref{fig:simul2} depicts the results of parametric and non-parametric statistical tests for simulated non-Gaussian-distributed data of $S=20$ subjects following either the fixed-effect model ({top} left panel) or the random-effects model ({top right} panel). For comparison, the results obtained on Gaussian-distributed data following a random-effects model are displayed in the {bottom} panel. Four different statistical tests are compared: 1) the random-effects inverse-variance-weighted sufficient-summary-statistic approach for the difference between class-conditional means, 2) the same test for the area under the non-parametric receiver-operating curve (AUC), 3) the naive summary-statistic approach in the form of a paired t-test between subject-wise means, and 4) its non-parametric equivalent, the Wilcoxon signed rank test. Note that for the naive summary-statistic approaches, the mean differences of each subject are treated as observations for a group-level paired t-test or Wilcoxon signed rank test, respectively. 

The figure shows that, as for Gaussian-distributed data, the inverse-variance-weighted sufficient-summary-statistic approach achieves considerably higher statistical power than the corresponding naive summary-statistic approaches. Furthermore, non-parametric approaches achieve a higher power for non-Gaussian-distributed data than their parametric equivalents assuming Gaussian-distributed data. This difference is particularly pronounced for the better performing inverse-variance-weighted sufficient-summary-statistic approaches. 
The difference for the naive summary approaches is much smaller, because subject-level averages tend to be more Gaussian according to the central limit theorem. In contrast, parametric approaches have only a very minor advantage over non-parametric ones for Gaussian-distributed data. Note further that, when the Gaussianity assumption of the parametric approaches is violated, spurious results can, in theory, not be ruled out. However, such effects are very small here.

\subsection{{Simulation 3: Correlation between subject-level effect size and variance}}
\label{sec:bias_example}
\begin{figure*}
\centering
\includegraphics[width =0.45\textwidth]{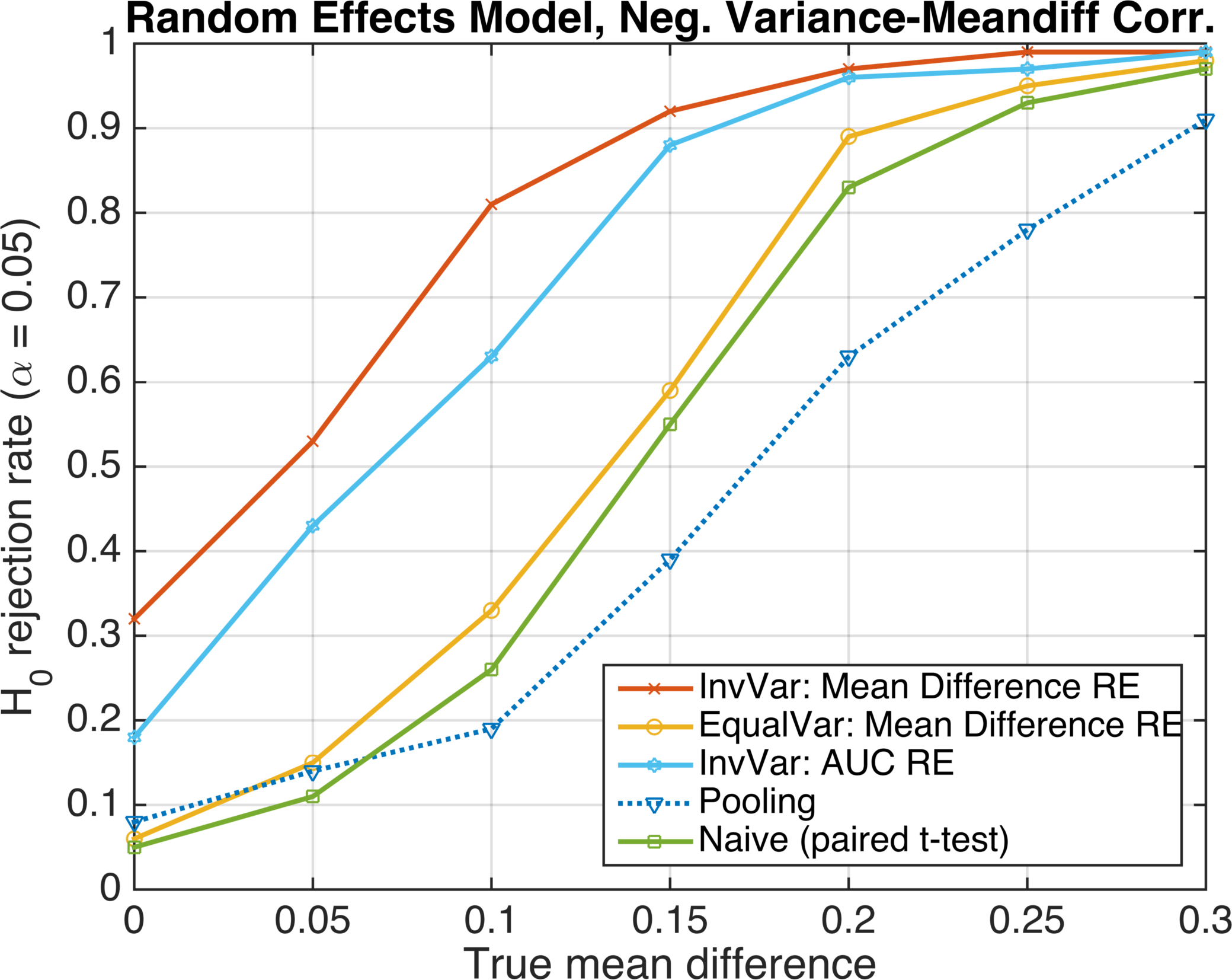} 
\includegraphics[width =0.45\textwidth]{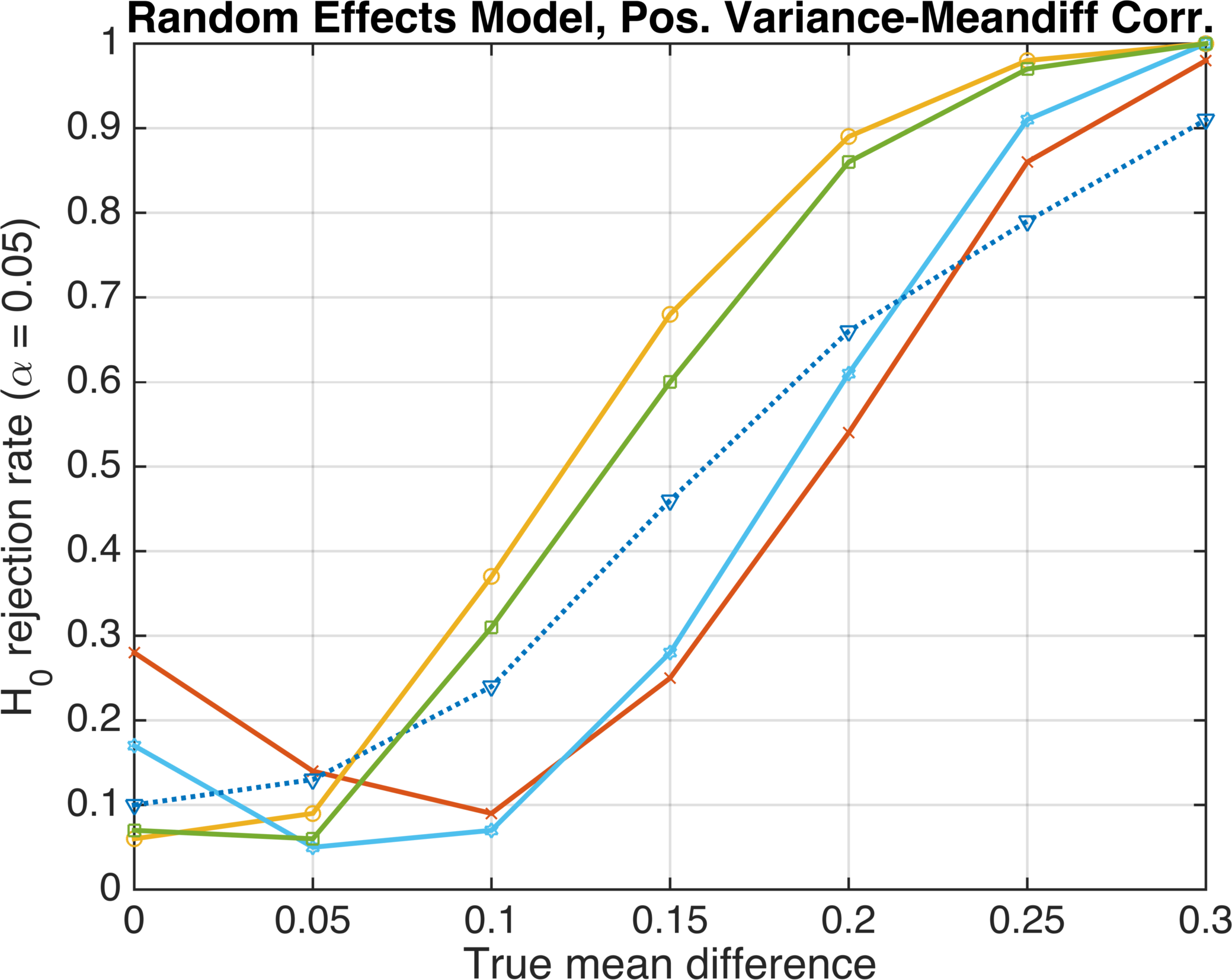} 
\caption{The probability of rejecting the null hypothesis $H_0: d = 0$ as a
function of the true mean difference $d$ of Gaussian-distributed simulated data from $S=20$ subjects. Left: data following a random-effects model and exhibiting a perfect negative correlation between the difference of the two group means and the standard deviation associated with this difference. Right: perfect positive correlation.}
\label{fig:simul3}
\end{figure*}

{In the presence of dependencies between subject-level effect sizes and corresponding variances, the resulting group-level effect size may become biased if inverse-variance weighting is used. To demonstrate this adverse effect, we simulated data exhibiting a perfect correlation between the difference of the two group means and the standard deviation associated with this difference, represented by the square root of Eq.~\eqref{eq:meanDifferenceVariance}. The between-subject variance was thereby kept at the same level as in the preceding random-effects simulations (cf. bottom-right panel of Figure~\ref{fig:simul1}, bottom panel of Figure~\ref{fig:simul2}, and corresponding texts). The left panel of Figure~\ref{fig:simul3} shows $H_0$-rejection rates for a negative correlation, implying that subjects with lower (negative) mean differences exhibit larger variability. In the inverse-variance-weighting approach, the influence of these subjects is down-weighted, which leads to an overestimation of the actual mean difference at the group-level. As a result, the number of false positive detections under the null hypothesis $H_0 : d = 0$ dramatically increases. The right panel of Figure~\ref{fig:simul3} shows analogous results for a positive correlation, implying that subjects with larger (positive) mean differences exhibit larger variability, and that inverse-variance weighting, consequently, introduces a negative bias on the group-level mean difference. Thus, when using a two-tailed statistical test, it may happen that a significant negative mean difference is found even in the presence of positive difference or the absence of any difference. This behavior is illustrated by the U-shape of the read and cyan curves. Note, however, that these problems do not occur if the sufficient-summary-statistics approach with equal weighting or the naive group-level t-test are used. }
  
\section{Analysis of emergency-braking-related brain activity}
\label{sec:realData}

We analyzed neuro- and myoelectrical activity of human participants during a simulated driving experiment. During the experiment, participants had the task to closely follow a computer-controlled lead vehicle. This lead vehicle would occasionally slow down abruptly, in which case the participant had to perform an emergency braking. The full study is described in \cite{haufe2011eeg}. Brain signals were acquired using 64 EEG electrodes (referenced to an electrode on the nose), while we here only report on the central EEG electrode Cz. Muscular activation of the lower right leg was acquired from two electromyographic (EMG) electrodes using a dipolar derivation. EEG and EMG Data were recorded from 18 participants in three blocks \`a 45 minutes. On average, clean data from 200 emergency situations were obtained from each participant (min: 123, max: 233). After filtering and sub-sampling to 100~Hz, the data were aligned (`epoched') relative to the onset of the braking of the lead vehicle as indicated by its brake light. For each time point relative to this stimulus, EEG and EMG measurements were contrasted with a sample of identical size that had been obtained from normal driving periods of each participant. {While for the present study only preprocessed and epoched data were used, original raw data are also publicly available\footnote{{\url{http://bnci-horizon-2020.eu/database/data-sets} (\#24)}}.}

\begin{figure}
\centering
\includegraphics[width =\linewidth]{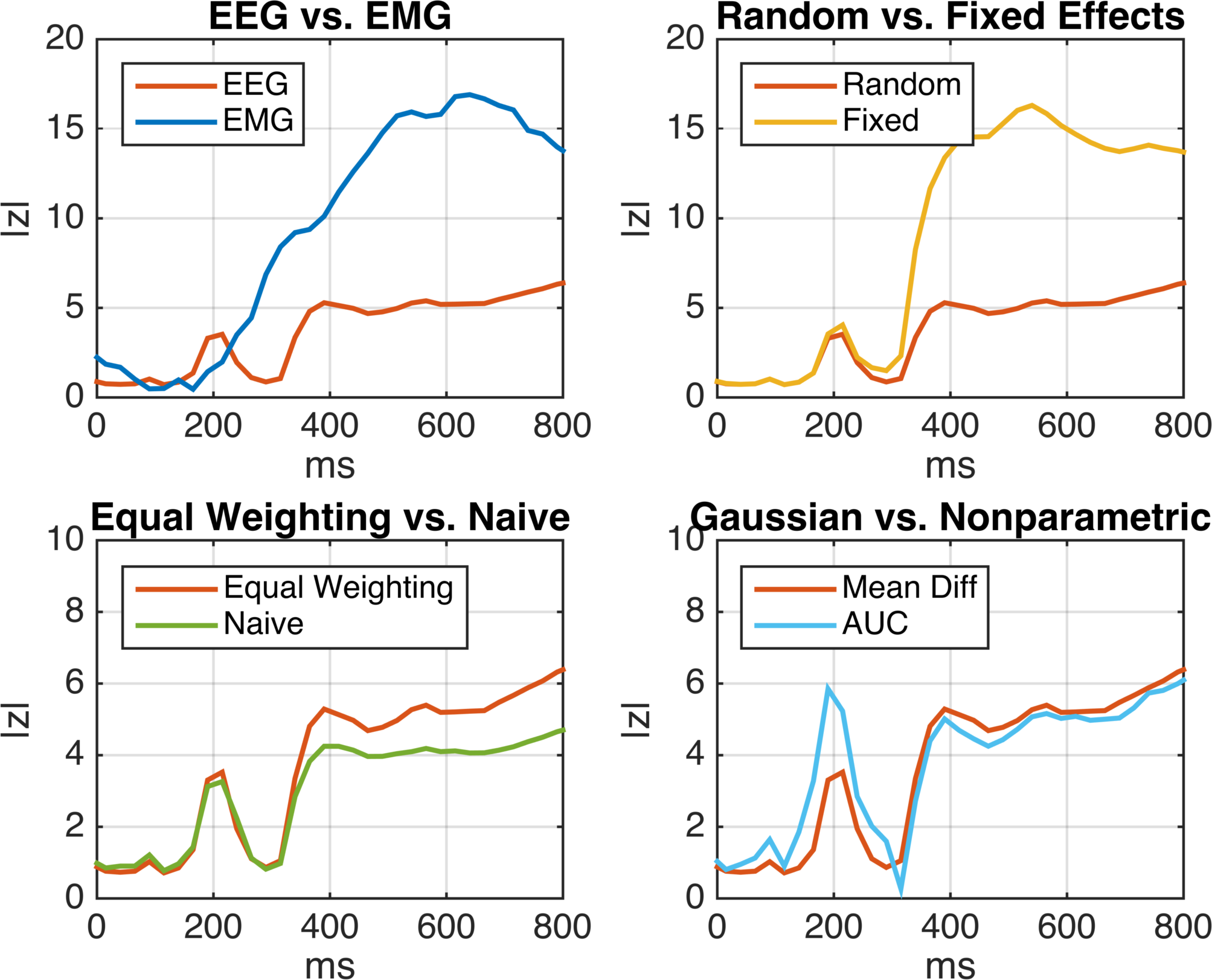} 
\caption{Analysis of event-related EEG (neural) and EMG (muscular) activity of $N=18$ car drivers during simulated emergency braking. Shown is the z-scaled difference between the mean activity during emergency braking situations and the mean activity during normal driving periods as a function of time after the emergency-initiating situation. Top Left: Comparison of EEG and EMG under the random-effects (RE) model using two-sample subject level t-tests and equal weighting. EEG displays a significant class separation at an earlier time than EMG, reflecting the logical order of the underlying perceptual decision-making process. Top Right: Comparison between the fixed-effect (FE) and RE models for EEG. The FE model displays inflated z-scores, indicating substantial but unaccounted between-subject variability. Bottom Left: Comparison of the naive summary statistic approach and the sufficient-summary-statistic approach using equal weighting for EEG (RE model). By taking the subject-level variances into account, the sufficient-summary-statistic approach achieves a clearer separation. Bottom Right: Comparison between the two-sample t-test and the non-parametric Wilcoxon-Mann-Whitney test for a group-level area under the ROC curve (AUC) greater than chance-level (RE model, equal weighting). Both tests lead to similar results, indicating that the distribution of samples is close to normal.
\label{fig:erp}}
\end{figure}

Figure~\ref{fig:erp} (top left) shows the deviation of EEG and EMG signals in emergency braking situations from signals obtained during normal driving periods as a function of time after stimulus. For each participant, the mean difference between the two driving conditions was computed (Eq.~\eqref{eq:meanDifference}). Assuming a random-effects model, the within-subject (i.e., within-participant) variance was estimated using Eq.~\eqref{eq:meanDifferenceVariance}, while the between-subject variance was estimated using Eq.~\eqref{eq:DersimonianLaird}. {We tested for Pearson correlations between subject-level mean differences and corresponding within-subject standard deviations and found strong significant positive correlations for almost all time points post-stimulus at both electrodes. As, under these circumstances, inverse-variance weighting is expected to produce biased results, we resorted to using the sufficient-summary-statistics approach in combination with equal weights for each subject.} Results are presented in terms of the absolute value of the group-level z-score, which was computed using equal weighting along the lines of Algorithm~\ref{algo:invVar}. It is apparent that the brain activity measured by EEG exhibits a significant amount of emergency-braking-related information at an earlier point in time than the activity measured at the right leg muscle, but is superseded in terms of class separability by the EMG later on. This result reflects the decision-making process that is taking place in the brain prior to the execution 
of the physical braking movement.

The top right panel of Figure~\ref{fig:erp} depicts the same EEG time course in comparison to the curve obtained under the fixed-effect model. Compared to the RE model, the FE model leads to an inflation of z-scores starting 300\,ms post-stimulus. Note that this is consistent with the result of Cochran's Q-test for effect size heterogeneity 
\citep{Cochran1954} indicating non-zero between-subject variability after 200\,ms post-stimulus ($p < 0.05$), but not before.

The bottom left panel of Figure~\ref{fig:erp} depicts the difference between the equal-weighting sufficient-summary-statistic approach and the naive summary-statistic approach implemented as a paired t-test for differences in the subject-wise means of the two conditions. As expected, the equal-weighting approach achieves a higher power than the naive approach (at least during later time points) by taking the subject-level variances into account. 

Finally, the bottom right panel of Figure~\ref{fig:erp} depicts the difference between subject-level two-sample t-tests and non-parametric AUC tests according to Eqs.~\eqref{eq:AU_transform} and \eqref{eq:StdAUC_general}. No substantial difference is found between the two except for a narrow time interval around 200~ms post-stimulus, in which the non-parametric test yields higher z-scores. Overall, this result suggests that the raw samples are approximately normal distributed, justifying the use of the parametric test.

\section{Discussion}
\label{sec:discussion}

In this paper we have provided a review of existing methods to assess the statistical significance of group-level effect sizes in data with nested structure. We demonstrated that simply pooling the data of all subjects is not a valid approach. The naive summary-statistic approach of performing a paired t-test on subject-level effect sizes is valid, but has suboptimal statistical power. With the {sufficient-summary-statistic} approach and Stouffer's method, we discussed two general strategies that combine the simplicity and low complexity of `naive' approaches with higher statistical power by using prior knowledge about the distributions and variances of the subject-level effect sizes. The benefit of these two strategies over the `naive' approaches was demonstrated in a set of simulations. Note that the degree of improvement due to using sufficient summary statistics depends on the number of trials per subject vs. the number of subjects. Therefore, differing observations can be found in the literature \citep[e.g.,][]{Beckmann2003,Mumford2009}.

The simulations as well as the presented real-data analysis also highlighted the necessity to account for between-subject variances through a random-effects analysis. A failure to do so results in underestimated p-values and the spurious detection of non-existing effects. Stouffer's method is a fixed-effects analysis, and thus provides a valid group-level test only if the assumption of zero between-subjects variance can be theoretically justified. In most practical cases, this is not the case \citep{Holmes1998,Field2003,stephan2009bayesian,Schmidt2009,allefeld2016valid}. We thus recommend the use of the {sufficient-summary-statistic} approach when the number of subjects is modest and the subject-wise variances can be reliably estimated. 

Importantly, while we here only considered data with two nesting levels, both Stouffer's method and the {sufficient-summary-statistic} approach naturally extend to hierarchies with arbitrary numbers of levels. For example, p-values derived from individual subjects of a study, e.g., using Stouffer's method, can again be combined at a higher level to test for consistent effects across multiple studies. In a similar way, group-level effects with variances derived from subject-level samples through Eq.~\eqref{eq:z_weightedmean} can be further combined into a higher-level average with known variance. 

\subsection{{Limitation of inverse-variance weighting}}
\label{sec:invvar_limitation}
{Our simulations demonstrated that inverse-variance weighting consistently outperformed all other approaches provided that no dependencies between subject-level effect sizes and their variances were present. Our analysis of real EEG data, however, also showed that such dependencies are not unlikely. In this example, participants with stronger emergency-related brain responses also showed larger variability. The opposite case is also conceivable, as participants with weaker responses (e.g. due to missing to process experimental stimuli) may exhibit more variable activity representing their less constrained mental state. In practice, it is, therefore, advisable to test for monotonous relationships between effect size and variance, for example using Spearman's rank correlation. If no such relationship is found, the inverse-variance-weighted sufficient-summary-statistics approach can be used. In the opposite case, we recommend the use of the sufficient-summary-statistics approach using equal weights, which still improves upon the naive group-level t-test. An alternative is to use weights proportional to the subject-level sample sizes \citep{marin2010weighting}.
Notwithstanding these considerations, dependencies between effect sizes and their variances must not be considered ubiquitous. By definition, they cannot occur in the presence of a fixed effect. Moreover, the variances of some effect size measures (e.g., Fisher-transformed Pearson correlations) only depend on the number of subject-level samples, and are constant if identical samples sizes are available for all subjects. In these settings, no systematic correlations and, for that matter, no biases can be expected, implying that inverse-variance weighting remains a valid and powerful approach.}

\subsection{Alternative definitions of fixed and random effects} 
The notions of `fixed' and `random' effects are used differently in different branches of statistics. See, for example, \cite{Gelman2005} for a discussion of five different definitions of `fixed' and `random' effects in the statistical literature. In ANOVA, the factor levels of a `random effect' are assumed to be randomly selected from a population, while the factor levels of a `fixed effect' are chosen by the experimenter. In contrast to the definition of a `fixed effect' used here (Eq.~\eqref{eq:FEModel}), the effect sizes of a `fixed effect' factor in ANOVA are allowed to differ across subjects.

Here we define a fixed effect (FE) to be constant across subjects, while a random effect (RE) is allowed to vary across subjects. The fundamental model underlying RE analysis is given by Eqs.~ \eqref{eq:REModel1} and \eqref{eq:REModel2}, while the FE model is defined in Eq.~\eqref{eq:FEModel}. These definitions are used in the meta-analysis literature \citep{Field2003,Borenstein2009,Card2011}, which 
contains most statistical discussion of between-subject variance estimators \citep{Dersimonian1986, Brockwell2001,Schmidt2009,Rukhin2013}. 

In parts of the neuroimaging literature, a different interpretation of the fixed-effect model is predominant \citep{Penny2006,Monti2011}. Here,
\begin{align}
\hat{\theta}_s = \theta_s + \epsilon_s,
\end{align}
where $\epsilon_s$ denotes the deviation of the subject's observed effect from the subject-specific true effect $\theta_s$, which is not modeled as a random variable. In this view, the subjects are not randomly drawn from a population, but are `fixed'. There is no overall population effect $\theta$ and the implicit null hypothesis behind the model is $H_0: 1/S \sum_{s=1}^S \theta_s = \theta_0$. This yields an alternative interpretation of the same analysis: a fixed-effect analysis allows one to draw valid inference on the mean effect -- but only for the specific mean of the observed subjects. Such an analysis would correspond to a case study, but a generalization to the population from which the subjects were drawn is not possible \citep{Penny2006}. In contrast, the fixed-effect model Eq.~\eqref{eq:FEModel} we assume throughout this paper allows such a generalization -- but the assumption of a constant effect across subjects has to be theoretically justified. 

\subsection{Nested multiple linear models} \label{sec:nestedlinear}

Another approach to handle nested data are nested linear models (also called hierarchical linear models, multi-level models or mixed linear models). These models extend the multiple linear regression model discussed in Section~\ref{sec:linear_model} to deal with nested data. Following \cite{hox2010multilevel}, this is done by introducing subject-specific regression coefficients $\beta_{k,s}, k \in \{0, \hdots, K\}, s \in \{1, \hdots, S\}$. The model for the $n$-th sample of subject $s$ then reads 
\begin{align}\label{eq:multlin}
y_{n, s} = \beta_{0, s} + x_{n,1} \beta_{1, s} + \hdots + x_{n,K} \beta_{K, s} + \epsilon_{n, s} \;.
\end{align}
The subject-specific coefficients are further expressed as 
\begin{align}\label{eq:multlin}
\beta_{k,s} = \gamma_{0, s} + \boldsymbol{\gamma}_{s}^\top \vec{z}_s + \varepsilon_{k, s} \;,
\end{align}
where $\gamma_{0, s}$ is a subject-specific intercept, $\vec{z}_s = (z_{1, s}, \hdots, z_{L, s})^\top$ models $L$ known subject-resolved independent variables $z_{l, s}$, $\boldsymbol{\gamma}_{s} = (\gamma_{1, s}, \hdots, \gamma_{L, s})^\top$ is a vector of corresponding coefficients modeling the influence of these variables on $\beta_{k,s}$, and $\varepsilon_{k, s}$ is group-level zero-mean noise. In this complete form, all coefficients are subject-specific. We therefore speak of a random-effects nested linear model. It is also conceivable that only some of the coefficients are subject-specific, while others are shared between subjects. For example, in some applications it may be reasonable to model subject-specific intercepts $\beta_{0, s}$, but identical effects $\beta_{k, 1} = \hdots = \beta_{k, S} = \beta_k$ for all subjects. A resulting model would be called a mixed-effects nested linear model. 

Nested linear models are very general and allow for more complex statistical analysis than the procedures for estimating and testing group-level effects discussed here. On the downside, the estimation of nested linear models is difficult because no closed-form solution exists in the likely case that the variances of the subject- and group-level noise terms are unknown \citep[e.g.,][]{chen2013linear}. Fitting a nested linear model using iterative methods is time consuming when the number of subjects and/or samples per subject is large, as all data of all subjects enter the same model. This is especially problematic when the number of models to be fitted is large, as, for example, in a mass-univariate fMRI context, where an individual model needs to be fitted for ten-thousands of brain voxels. 

When only the group-level effect is of interest, the presented sufficient-summary-statistic approach is the more practical and computationally favorable alternative. In this approach, regression coefficients $\hat{\beta}_{k,s}$ are estimated at the subject level, which bears the advantage that the global optimum for each subject can be found analytically in a computationally efficient manner. As the individual $\hat{\beta}_{k,s}$ are normal distributed with variance given in Eq.~\eqref{eq:beta_variance}, they can then be combined, e.g., using the inverse-variance-weighting scheme. This approach is mathematically equivalent to a nested-linear model analysis when the covariances are known \citep{Beckmann2003}.  For these reasons, we here refrained from a deeper discussion of nested linear models. The interested reader is referred to, for example, \cite{Quene2004,Woltman2012,hox2010multilevel,chen2013linear}.

\subsection{Resampling and surrogate-data approaches}
While the variances of the effect size measures discussed here can be derived analytically, this may not be the case in general. However, given sufficient data, the variance of the observed effect $\hat{\theta}$ can always be estimated through resampling procedures such as the bootstrap or the jackknife \citep{efron1982jackknife}. Assuming an approximately normal distribution for $\hat{\theta}$, the inverse-variance-weighting approach can be applied. 

For some types of data such as time series, the subject-level i.i.d. assumption underlying most statistical procedures discussed here is, however, violated. For such dependent samples, the variance of an observed effect $\hat{\theta}$ -- be it analytically derived or obtained through a resampling procedure under the i.i.d. assumption -- is underestimated. This problem can be addressed through sophisticated resampling techniques which accommodate dependent data structure. A detailed describtion of these techniques can be found, for example, in \cite{Kreiss2011,Lahiri2013}.

For some types of analysis questions, it is not straightforward to determine the expected effect under the null hypothesis $\theta_0$. A potential remedy to this problem is the method of surrogate data.  
Surrogate data are artificial data that are generated by manipulating the original data in a way such that all crucial properties (including the dependency structure of the samples) are maintained except for the effect that is measured by $\theta$. As such, surrogate data can provide an empirical distribution for $\hat{\theta}$ under the null hypothesis. This may be used to derive subject-level p-values, which can be subjected to Stouffer's method to test for population effects under the fixed-effect model. Originally introduced in the context of identifying nonlinearity in time series \citep{theilersurro}, variants of this approach are increasingly often applied to test for interactions between neural time series \citep[e.g.,][]{honey2012slow,haufe2016simulation,Haufe_Elucidating2017}.

\subsection{Multivariate statistics} 
In the present paper we assumed that only a single effect is measured for each subject. However, massively multivariate data are common especially in neuroimaging, where brain activity is typically measured at hundreds or even thousands of locations in parallel. When (group) statistical inference is performed jointly for multiple measurement channels, the resulting group-level p-values must be corrected for multiple comparisons {using, e.g., methods described in \cite{Genovese2002,Nichols2003,Pernet2015}.}

Another way to perform inference for multivariate data is to use inherently multivariate effect size measures such as canonical correlations, coefficients of multivariate linear models, {the accuracy of a classifier \citep[e.g.,][]{Haxby2001,Norman2006},} or more generally univariate effect size measures that are calculated on optimal linear combination of the measurement channels \citep[e.g.,][]{DaeBieSamHauGolGunVilFazMue2015,HauMeiGoeDaeHayBlaBie13}. However, most multivariate statistics involve some sort of model fitting. If the number of data channels is high compared to the number of samples, overfitting may occur, and may bias the expected value of the effect under the null hypothesis. One way to avoid that bias by splitting the data into training and test parts, where the training set is used to fit the parameters of the multivariate model, while the actual statistical test is carried out on the test data using the predetermined model parameters \citep{Lemm2011387}. 



\subsection{Activation- vs. information-like effect size measures.}
A distinction is made in the neuroimaging literature between so-called `activation-like' and `information-like' effect size measures. 
\cite{allefeld2016valid} argue that measures that quantify the presence of an effect without a notion of directionality (that is, are `information-like') cannot be subjected to a subsequent random-effects group-level analysis, because their domain is bounded from below by what would be expected under the null hypothesis of no effect. Their arguments refers in particular to the practice of plugging single-subject classification accuracies into a group-level paired t-test. Because the true single-subject classification accuracies can never be below chance-level, the group-level null hypothesis being tested is the fixed-effect hypothesis of no effect in any subject. {Another problem with`information-like' measures is that certain confounds are not appropriately controlled for, because confounding effects of different direction do not cancel each other out \citep{Todd2013}.} For the current investigation, {these issues are}, however, of minor importance, as, except for the coefficient of determination, all effect size measures discussed here are directional and therefore `activation-like'.  



\section{Conclusion}
\label{sec:conclusions}
In this paper, we have reviewed practical approaches to conduct statistical inference on group-level effects in nested data settings, and have demonstrated their properties on simulated and real neuroimaging data. {With the sufficient-summary-statistic approach, we highlighted an approach that combines computational simplicity with favorable statistical properties.} We have furthermore provided a practical guideline for using this approach in conjunction with some of the most popular measures of statistical effects.


%



\section*{Statement regarding ethical research practices}
This work includes human EEG data from a previously published study \citep{haufe2011eeg}. For this study, healthy adult volunteers from the general public were monitored with non-invasive EEG while using a driving simulator. Participants received monetary compensation proportional to the time spent in the laboratory, and signed an informed consent sheet stating that they could withdraw at any time from the experiment without negative consequence. The experiments were conducted in accordance with ethical principles
that have their origin in the Declaration of Helsinki. The protocol was, however, not sent to an institutional review board (IRB) for approval. This was due to the absence of IRBs dealing with non-medical studies in Berlin at that time.

\section*{Conflict of Interest Statement}

The authors declare that the research was conducted in the absence of any commercial or financial relationships that could be construed as a potential conflict of interest.

\section*{Author Contributions}

Both authors contributed substantially to all aspects of the work including the conception and design of the work, the acquisition, analysis, and interpretation of data for the work, drafting of the work and revising it critically for important intellectual content. Both authors approved the final version to be published and agree to be accountable for the content of the work.

\section*{Funding}
SH was supported by a Marie Curie Individual International Outgoing Fellowship (grant No. PIOF-GA-2013-625991) within the 7th European Community Framework Programme.

\section*{Acknowledgments}
We would like to thank the Reviewers for providing thoughtful comments that helped to improve the paper.

\bibliographystyle{frontiersinSCNS_ENG_HUMS} 
\bibliography{groupstat}



\end{document}